\input eplain
\input BMmacs

%
%
%
%
%
%
%
%
%
\def\YEAR{\year}\newcount\VOL\VOL=\YEAR\advance\VOL by-1995
\def\firstpage{1}\def\lastpage{1000}
\def\received{}\def\revised{}
\def\communicated{}
\nopagenumbers

\hoffset=1.91cm\voffset=2.7cm 

\hsize=12.0cm\vsize=19.0cm

\font\eightrm=cmr8
\font\caps=cmcsc10                    
\font\Caps=cmcsc10 scaled \magstep1   

%
\pageno=\firstpage\def\folio{\rm\number\pageno}\def\headfoot{25pt}
\def\DocMath{}
\def\makeheadline{
    \vbox to 0pt{\vskip-\headfoot\line{\vbox to8.5pt{}\the\headline}\vss}
    \nointerlineskip}
\headline={\ifnum\pageno=\firstpage{\DocMath\hfill\llap{\folio}}%
           \else{\ifodd\pageno\rightheadline\else\leftheadline\fi}\fi}
\def\rightheadline{\hfill {\caps\rightheadtext}    \hfill \llap{\folio}}
\def\leftheadline {\rlap{\folio} \hfill {\caps\leftheadtext} \hfill }
\def\leftheadtext{\ifnum\pageno>\lastpage\else\SAuthor\fi}
\def\rightheadtext{\STitle}
\def\TSkip{\bigskip}
\newbox\TheTitle{\obeylines\gdef\GetTitle #1
\ShortTitle  #2
\SubTitle    #3
\Author      #4
\ShortAuthor #5
\EndTitle
{\setbox\TheTitle=\vbox{\baselineskip=20pt\let\par=\cr\obeylines%
\halign{\centerline{\Caps##}\cr\noalign{\medskip}\cr#1\cr}}%
	\copy\TheTitle\TSkip\TSkip%
\def\next{#2}\ifx\next\empty\gdef\STitle{#1}\else\gdef\STitle{#2}\fi%
\def\next{#3}\ifx\next\empty%
    \else\setbox\TheTitle=\vbox{\baselineskip=20pt\let\par=\cr\obeylines%
    \halign{\centerline{\caps##} #3\cr}}\copy\TheTitle\TSkip\TSkip\fi%
\centerline{\caps #4}\TSkip\TSkip%
\def\next{#5}\ifx\next\empty\gdef\SAuthor{#4}\else\gdef\SAuthor{#5}\fi%
\ifx\received\empty\relax
    \else\centerline{\eightrm Received: \received}\fi%
\ifx\revised\empty\TSkip%
    \else\centerline{\eightrm Revised: \revised}\TSkip\fi%
\ifx\communicated\empty\relax
    \else\centerline{\eightrm Communicated by \communicated}\fi\TSkip\TSkip%
\catcode'015=5}}\def\Title{\obeylines\GetTitle}
\def\Abstract{\begingroup\narrower
    \parskip=\medskipamount\parindent=0pt{\caps Abstract. }}
\def\EndAbstract{\par\endgroup\TSkip}

\long\def\MSC#1\EndMSC{\def\arg{#1}\ifx\arg\empty\relax\else
     {\par\narrower\noindent%
     2000 Mathematics Subject Classification: #1\par}\fi}

\long\def\KEY#1\EndKEY{\def\arg{#1}\ifx\arg\empty\relax\else
	{\par\narrower\noindent Keywords and Phrases: #1\par}\fi\TSkip}

\newbox\TheAdd\def\Addresses{\vfill\copy\TheAdd\vfill
    \ifodd\number\lastpage\vfill\eject\phantom{.}\vfill\eject\fi}
{\obeylines\gdef\GetAddress #1
\Address #2 
\Address #3
\Address #4
\EndAddress
{\def\xs{4.3truecm}\parindent=0pt
\setbox0=\vtop{{\obeylines\hsize=\xs#1\par}}\def\next{#2}
\ifx\next\empty 
     \setbox\TheAdd=\hbox to\hsize{\hfill\copy0\hfill}
\else\setbox1=\vtop{{\obeylines\hsize=\xs#2\par}}\def\next{#3}
\ifx\next\empty 
     \setbox\TheAdd=\hbox to\hsize{\hfill\copy0\hfill\copy1\hfill}
\else\setbox2=\vtop{{\obeylines\hsize=\xs#3\par}}\def\next{#4}
\ifx\next\empty\ 
     \setbox\TheAdd=\vtop{\hbox to\hsize{\hfill\copy0\hfill\copy1\hfill}
                \vskip20pt\hbox to\hsize{\hfill\copy2\hfill}}
\else\setbox3=\vtop{{\obeylines\hsize=\xs#4\par}}
     \setbox\TheAdd=\vtop{\hbox to\hsize{\hfill\copy0\hfill\copy1\hfill}
	        \vskip20pt\hbox to\hsize{\hfill\copy2\hfill\copy3\hfill}}
\fi\fi\fi\catcode'015=5}}\gdef\Address{\obeylines\GetAddress}

\hfuzz=0.1pt\tolerance=2000\emergencystretch=20pt\overfullrule=5pt
\Title

Special subvarieties arising from
families of cyclic covers of the projective line

\ShortTitle 

Special subvarieties and families of cyclic covers

\SubTitle   

{\it Dedicated to Frans Oort on the occasion of his 75th birthday}

\Author 

Ben Moonen

\ShortAuthor 
\EndTitle
\Abstract 

We consider families of cyclic covers of $\mP^1$, where we fix the covering group and the local monodromies and we vary the branch points. We prove that there are precisely twenty such families that give rise to a special subvariety in the moduli space of abelian varieties. Our proof uses techniques in mixed characteristics due to Dwork and Ogus.

\EndAbstract
\MSC 

11G15, 14H40, 14G35

\EndMSC
\KEY 

Special subvarieties, Jacobians, complex multiplication

\EndKEY
\Address 
\Address
University of Amsterdam
KdV Institute for Mathematics
PO Box 94248
1090 GE Amsterdam
The Netherlands
\Address
\Address
\EndAddress


\intro

\introssect
According to a conjecture of Coleman, if we fix a genus $g \geq 4$ and consider smooth projective curves~$C$ over~$\mC$ of genus~$g$, there should be only finitely many such curves, up to isomorphism, such that $\Jac(C)$ is an abelian variety of CM type. This conjecture is known to be false for $g \in \{4,5,6,7\}$, by virtue of the fact that for these genera there exist special subvarieties $S \subset \cA_g$ (also known as subvarieties of Hodge type) of positive dimension that are contained in the (closed) Torelli locus and that meet the open Torelli locus.

All known examples of such special subvarieties~$S$ arise from families of cyclic covers of~$\mP^1$. As input for this we fix integers $m \geq 2$ and $N \geq 4$, together with an $N$-tuple $a = (a_1,\ldots,a_N)$; then we consider cyclic covers $C_t \to \mP^1$ with covering group $\mZ/m\mZ$, branch points $t_1,\ldots,t_N$ in~$\mP^1$, and with local monodromy~$a_i$ about~$t_i$. Varying the branch points we obtain an $(N-3)$-dimensional closed subvariety $Z = Z(m,N,a) \subset \cA_g$ given by the Jacobians of the curves~$C_t$, and for certain choices of $(m,N,a)$ it can be shown that $Z$ is a special subvariety. 

The main purpose of this paper is to give a list of monodromy data $(m,N,a)$ such that $Z(m,N,a)$ is special, and to prove that there are no further such examples.

\introssect
By construction, the Jacobians $J_t$ come equipped with an action of the group ring $\mZ[\mu_m]$. This action defines a special subvariety $S(\mu_m) \subset \cA_g$ that contains~$Z$ and whose dimension can be easily calculated in terms of the given monodromy data $(m,N,a)$. We always have $N-3 \leq \dim S(\mu_m)$, and if equality holds $Z = S(\mu_m)$ is special. To search for triples $(m,N,a)$ for which $\dim S(\mu_m) = N-3$ is something that can be done on a computer, and it is therefore surprising that this appears not to have been done until recently. Thus, while examples in genera~$4$ and~$6$ have been given by de Jong and Noot in~[\ref{dJN}]---and in fact, these examples at least go back to Shimura's paper~[\ref{Shimura}]---examples with $g=5$ and $g=7$ were found only much later; see Rohde~[\ref{Rohde}]. At any rate, with the help of a computer program we find twenty examples with $\dim S(\mu_m) = N-3$. The main result of this paper is that this list is complete:

\introssect
{\it Main theorem. --- Consider monodromy data $(m,N,a)$ as above. Then the closed subvariety $Z(m,N,a) \subset \cA_{g,\mC}$ is special if and only if $(m,N,a)$ is equivalent to one of the twenty triples listed in\/~{\rm Table~\ref{tablespecials}}.}

\introssect
Note that if $N-3 < \dim S(\mu_m)$, which means that the inclusion $Z \subset S(\mu_m)$ is strict, $Z$ could a priori still be a special subvariety. The main point of our result is that even in these cases we are able to prove that $Z$ is not special.

\introssect
To get a feeling for the difficulty of the problem, consider, as an example, the family of (smooth projective) curves of genus~$8$ given by $y^{10} = x(x-1)(x-t)^2$, where $t \in T = \mP^1_\mC\setminus\{0,1,\infty\}$ is a parameter. The corresponding family of Jacobians $J \to T$ decomposes, up to isogeny, as a product $J^\old \times J^\new$ of two abelian fourfolds. (The old part comes from the quotient family $u^5 = x(x-1)(x-t)^2$; the new part is the family of Pryms.) Both $J^\old$ and~$J^\new$ give rise to a $1$-dimensional special subvariety in~$\cA_4$, and they both admit an action by $\mZ[\zeta_5] = \mZ[\zeta_{10}]$, with the same multiplicities on the tangent spaces at the origin. A priori it might be true that $J^\old$ and~$J^\new$ are isogenous, in which case the family $J \to T$ would define a special subvariety in~$\cA_8$. Our theorem says, in this particular example, that this does not happen.

{}From the perspective of special points and the Andr\'e-Oort conjecture, one could state the problem, in the example at hand, as follows. There are infinitely many values for $t$ such that $J^\old_t$ is of CM type; likewise for~$J^\new_t$. Do there exist infinitely many~$t$ such that $J^\old_t$ and~$J^\new_t$ are simultaneously of CM type? Our theorem, combined with a result of Yafaev~[\ref{Yafaev}], implies the answer is {\it no\/}, at least if we assume the Generalized Riemann Hypothesis. See Corollary~\ref{CorMain}.

\introssect
Our proof of the main theorem uses techniques from characteristic~$p$ and is based on a method that in some particular cases has already been used by de Jong and Noot in [\ref{dJN}], Section~5. For a suitable choice of a prime number~$p$, the Jacobians in our family will, at least generically, have ordinary reduction in characteristic~$p$. Moreover, if we assume $Z(m,N,a) \subset \cA_g$ to be special then by a result of Noot~[\ref{NootThesis}] we can arrange things in such a way that the canonical liftings of these ordinary reductions are again Jacobians. Already over the Witt vectors of length~$2$ this is a very restrictive condition that, using the results of Dwork and Ogus in~[\ref{DworkOgus}], can be turned into something computable.

In individual cases (such as the example sketched above, or the examples treated in~[\ref{dJN}]), these techniques give a rather effective method to prove that some given family of curves does not give rise to a special subvariety in~$\cA_{g,\mC}$. To do this for arbitrary data $(m,N,a)$ is much harder, and some perseverance is required to deal with the combinatorics that is involved. It would be interesting to have a purely Hodge-theoretic proof of the main theorem; as it is, we do not even have a Hodge-theoretic method that works well in individual examples.

\introssect
Let us give an overview the contents of the individual sections. In Section~\ref{PrelimsSpecial} we quickly review the notion of a special subvariety and we summarize some facts we need. In Section~\ref{Setup} we introduce the families of curves that we want to consider. In Section~\ref{SpecialSubvar} we discuss the special subvariety $S(\mu_m)$ that contains $Z(m,N,a)$, we give the list of data $(m,N,a)$ for which $Z(m,N,a) = S(\mu_m)$ and we state the main results. Sections~\ref{DecompVHS} and~\ref{canlifts} contain the main technical tools for the proof. In Section~\ref{DecompVHS} we discuss how the VHS associated with our family of curves decomposes and we give some results about the monodromy of the summands. In Section~\ref{canlifts} we briefly review the techniques of~[\ref{DworkOgus}] and we give some refinements that we need in order to deal with hyperelliptic families. In Sections~\ref{PfN=4} and~\ref{PfN>4}, finally, we prove the main theorem, first in the case $N=4$, then for $N \geq 5$. We refer to the beginning of Section~\ref{PfN=4} for a brief explanation of how the proof works. 

\introssect
{\it Acknowledgements.\/} I thank Frans Oort for his encouragement and for stimulating discussions. It is a great pleasure to dedicate this paper to him. I thank Dion Gijswijt and Fokko van de Bult for helpful suggestions related to the combinatorics involved in the proof of Lemma~\ref{newpartnonrig}. Further I thank Maarten Hoeve, who wrote two Python programs (available on the author's webpage) to test some of the properties we need. The results of one of these programs are used in the proof of Lemma~\ref{newpartnonrig}. Finally, I thank the referee for his or her comments on the paper and for some very helpful suggestions regarding the exposition.

\introssect
{\it Notation.\/} If $x$ is a real number, $\frpart{x} = x-\lfloor x\rfloor$ denotes its fractional part. For $a \in \mZ/m\mZ$ or $a \in \mZ$ we denote by $\repr{a}_m$ the unique representative of the class~$a$ (resp.\ the class $a \bmod m$) in $\{0,\ldots,m-1\}$.

\section{Some preliminaries about special subvarieties}{PrelimsSpecial}
\bigskip

\noindent
In this section we work over~$\mC$. We recall the notion of a special subvariety. For further details we refer to~[\ref{BMLin1}], where the terminology ``subvarieties of Hodge type'' was used, and~[\ref{NootBourb}].

\subsection
\label{SiegelMod}
Consider the moduli space $\cA_{g,[n]}$ of $g$-dimensional principally polarized abelian varieties with a level~$n$ structure, for some $n \geq 3$. We first briefly recall the description of this moduli space as a Shimura variety.

Let $V_\mZ := \mZ^{2g} \subset V := \mQ^{2g}$, and let $\Psi \colon V_\mZ \times V_\mZ \to \mZ$ be the standard symplectic form. Let $G := \CSp(V_\mZ,\Psi)$ be the group of symplectic similitudes. Let $\mS := \Res_{\mC/\mR}\, \mG_m$ be the Deligne torus, and let $\gH$ be the space of homomorphisms $h \colon \mS \to G_\mR$ that define a Hodge structure of type $(-1,0) + (0,-1)$ on~$V_\mZ$ for which $\pm (2\pi i) \cdot \Psi$ is a polarization. The pair $(G_\mQ,\gH)$ is a Shimura datum, and $\cA_g$ can be described as the associated Shimura variety. Concretely, if $K_n := \bigl\{g \in G(\hat{\mZ}) \bigm| g \equiv 1 \bmod n\bigr\}$ then $\cA_{g,[n]}(\mC) \cong G(\mQ)\backslash \gH \times G(\Af)/K_n$.

In order to define special subvarieties, consider an algebraic subgroup $H \subset G_\mQ$ such that
$$
Y_H := \bigl\{h \in \gH \bigm| \hbox{$h$ factors through $H_\mR$}\bigr\}
$$
is non-empty. The group $H(\mR)$ acts on $Y_H$ by conjugation. It can be shown (see [\ref{BMthesis}], Section~I.3, or [\ref{BMLin1}], 2.4) that $Y_H$ is a finite union of orbits under $H(\mR)$. We remark that the condition that $Y_H$ is non-empty imposes strong restrictions on~$H$; it implies, for instance, that $H$ is reductive. If $Y^+ \subset Y_H$ is a connected component and $\eta K_n \in G(\Af)/K_n$, the image of $Y^+ \times \{\eta K_n\}$ in $\cA_{g,[n]}$ is an algebraic subvariety. 

\subsection
{\it Definition.\/} --- A closed irreducible algebraic subvariety $S \subset \cA_{g,[n]}$ is called a {\it special subvariety\/} if there exist $Y^+ \subset Y_H$ and $\eta K_n \in G(\Af)/K_n$ as above such that $S(\mC)$ is the image of $Y^+ \times \{\eta K_n\}$ in $\cA_{g,[n]}(\mC)$. 
\medskip

We refer to [\ref{BMLin1}] and~[\ref{NootBourb}] for alternative descriptions and basic properties of special subvarieties.

\subsection
\label{specstacks}
As level structures play no role in what we are doing, we prefer to state our results in terms of the moduli stack~$\cA_g$. Of course this is a stack, not a variety, but we allow ourselves to abuse terminology and speak of special subvarieties in~$\cA_g$ (over~$\mC$), where ``special substack'' would perhaps be more correct. By definition, then, a special subvariety $S \subset \cA_g$ is a closed, reduced and irreducible algebraic substack such that for some (equivalently: any) $n \geq 3$ the irreducible components of the inverse image of~$S$ under the natural map $\cA_{g,[n]} \to \cA_g$ are special. (It is in fact enough to require that some irreducible component is special.)

\subsection
\label{IsogenyFactors}
Let $f \colon X \to T$ be a family of $g$-dimensional abelian varieties over an irreducible non-singular complex algebraic variety~$T$. Choose a principal polarization~$\lambda$ of~$X/T$, and let $\phi \colon T \to \cA_g$ be the resulting morphism. Whether or not the closure of~$\phi(T)$ is a special subvariety of~$\cA_g$ only depends on the isogeny class of the generic fibre of~$X/T$ as an abelian variety over the function field of~$T$. If we replace~$X/T$ by an isogenous family or choose a different~$\lambda$, the morphism~$\phi$ is replaced by some other morphism $\phi^\prime \colon T \to \cA_g$, but if $\ol{\phi(T)}$ is special then so is $\ol{\phi^\prime(T)}$. (This is even true without the assumption that $\lambda$ is a principal polarization, but we shall not need this.) Without ambiguity we may therefore say that $X/T$ is special if $\ol{\phi(T)} \subset \cA_g$ is special for some choice of a polarization.

One property we shall use is that if $X/T$ is isogenous to $Y_1 \times Y_2$ and $X$ is special, the factors~$Y_i$ are both special. The converse is not true: if $Y_1/T$ and~$Y_2/T$ are special, it is not necessarily the case that $(Y_1 \times Y_2)/T$ is special.

\subsection
\label{MTfactors}
Let $f \colon X \to T$ be as in~\ref{IsogenyFactors}. Consider the $\mQ$-VHS over~$T$ with fibres the first cohomology groups $H^1(X_t,\mQ)$. Let $b \in T(\mC)$ be a Hodge-generic point for this VHS, and let $M \subset \GL\bigl(H^1(X_b,\mQ)\bigr)$ be the generic Mumford-Tate group of the family. Choose a principal polarization~$\lambda$ of~$X/T$, let $\phi \colon T \to \cA_g$ be the resulting morphism, and write $Z := \ol{\phi(T)}$. Finally, let $S \subset \cA_g$ be the smallest special subvariety containing~$Z$.

In order to calculate the dimension of~$S$, it suffices to know the adjoint real group~$M_\mR^\ad$. More precisely, if $M_\mR^\ad = Q_1 \times \cdots \times Q_r$ is the decomposition of this group as a product of simple factors, $\dim(S) = \sum_{i=1}^r d(Q_i)$, where $d(Q_i)$ is a contribution that only depends on the isomorphism class of the simple group~$Q_i$. The only cases that are relevant for us are that\par
\itemize
\item{$\bullet$} $d(Q) = 0$ if $Q$ is anisotropic (compact);
\item{$\bullet$} $d(Q) = pq$ if $Q \cong \PSU(p,q)$;
\item{$\bullet$} $d(Q) = h(h+1)/2$ if $Q \cong \PSp_{2h}$.
\stopitemize

\section{The setup}{Setup}
\bigskip

\noindent
In this section, given data $(m,N,a)$ as in the introduction (see~\ref{CoverTsetup} below), we construct a family of cyclic covers of $\mP^1$ over some base scheme~$T$. For later purposes we shall do this over a ring~$R$ of finite type over~$\mZ$.

\subsection
\label{CoverTsetup}
Let $m$ and $N$ be integers with $m \geq 2$ and $N \geq 2$, and consider an $N$-tuple of positive integers $a = (a_1,\ldots,a_N)$ such that $\gcd(m,a_1,\ldots,a_N) = 1$. We further require that $a_i \not\equiv 0$ modulo~$m$ for all~$i$ and $\sum_{i=1}^N a_i \equiv 0$ modulo~$m$. The triple $(m,N,a)$ serves as input for our constructions.

We call two such triples $(m,N,a)$ and $(m^\prime,N^\prime,a^\prime)$ equivalent if $m=m^\prime$ and $N=N^\prime$ and if the classes of $a$ and~$a^\prime$ in $(\mZ/m\mZ)^N$ are in the same orbit under $(\mZ/m\mZ)^* \times \gS_N$. Here we let $(\mZ/m\mZ)^*$ act diagonally by multiplication, and the symmetric group~$\gS_N$ acts by permutation of the indices.

In what follows we shall usually assume $N \geq 4$ but for some arguments it is useful to allow $N$ to be $2$ or~$3$.

\subsection
\label{FamilyConstr}
Let $(m,N,a)$ be a triple as in~\ref{CoverTsetup}. Let $R$ be the ring $\mZ[1/m,u]/\Phi_m$ with $\Phi_m$ the $m$th cyclotomic polynomial. We write $\zeta \in R$ for the class of~$u$; it is a root of unity of order~$m$. We embed $R$ into~$\mC$ by sending~$\zeta$ to $\exp(2\pi i/m)$. The element~$\zeta$ defines an isomorphism of $R$-group schemes $(\mZ/m\mZ)_R \isomarrow \mu_{m,R}$ by $(b \bmod m) \mapsto \zeta^b$.

Let $U \subset (\mA_R^1)^N$ be the complement of the big diagonals. In other words, $U$ is the $R$-scheme of ordered $N$-tuples $(t_1,\ldots,t_N)$ of distinct points in~$\mA^1$. Let $B \subset \mP^2_U$ be the projective curve over~$U$ obtained as the Zariski closure of the affine curve whose fibre over a point $(t_1,\ldots,t_N)$ is given by
$$
y^m = (x-t_1)^{a_1} \cdots (x-t_N)^{a_N} = \prod_{i=1}^N\, (x-t_i)^{a_i}\, .
$$ 
We have a $\mu_m$-action on $B$ over~$U$ by $\zeta \cdot (x,y) = (x,\zeta \cdot y)$. The rational function~$x$ defines a morphism $\pi_B\colon B \to \mP^1_U$.

There exist an open subscheme $T \subset U$, a smooth proper curve $f \colon C \to T$ equipped with an action of $\mu_{m,T}$, and a $\mu_{m,T}$-equivariant morphism $\rho \colon C \to B_T$, such that for every point $t \in T$ the morphism on fibres $\rho_t \colon C_t \to B_t$ is a normalization of~$B_t$. Let $\pi := \pi_B \circ \rho \colon C \to \mP^1_T$, which is a finite morphism that realizes $\mP^1_T$ as the quotient of~$C$ by the action of~$\mu_{m,T}$. If the context requires it, we include the data $(m,N,a)$ in the notation, writing $C = C(m,N,a)$ for instance. Further we write $J \to T$ for the Jacobian of $C$ over~$T$.

If $k$ is a field and $t = (t_1,\ldots,t_N) \in T(k)$, we can also describe $\pi_t\colon C_t \to \mP^1_k$ as the $\mu_m$-cover of~$\mP^1_k$ with branch points $t_1,\ldots,t_N$ and local monodromy about~$t_i$ given by the element $\zeta^{a_i} \in \mu_m$. 

The assumption that $\gcd(m,a_1,\ldots,a_N) =1$ implies that the fibres of $f\colon C\to T$ are geometrically irreducible. Let $r_i := \gcd(m,a_i)$. The Hurwitz formula gives that the fibres have genus
$$
g = 1 + {(N-2)m - \sum_{i=1}^N r_i \over 2}\, ,
\eqlabel{gformula}
$$
so we obtain a morphism $\psi \colon T \to \cM_g$ over~$R$. Define 
$$
\phi \colon T \to \cA_g
\eqlabel{phiDef}
$$ 
to be the composition of~$\psi$ with the Torelli morphism. Up to isomorphism, the morphism~$\psi$, and hence also~$\phi$, only depends on the equivalence class of the triple $(m,N,a)$ for the equivalence relation defined in~\ref{CoverTsetup}.

\subsection
{\it Remark.\/} --- For computational purposes we have restricted our attention to the case where all branch points~$t_i$ are in $\mA^1 \subset \mP^1$, and we have not fixed any of these branch points, thereby creating some redundancy and excluding families where one of the branch points is fixed to be the point at~$\infty$. For our main results, the only thing that really matters is the closure of the image of $\phi_\mC \colon T_\mC \to \cA_{g,\mC}$. For instance, while the family of curves $y^{10} = x(x-1)(x-\lambda)$ is not among the families we consider (the point at~$\infty$ being a branch point), the subvariety of~$\cA_9$ we obtain from this family is the same as the one obtained by taking $m=10$, with $N=4$ and $a=(1,1,1,7)$.

\subsection
{\it Convention.\/} --- In what follows we shall in several steps replace the base scheme~$T$ by a subscheme. In such a case, it will be understood that we replace~$C$ by its restriction to the new base scheme, and we again write $f \colon C \to T$ for the curve thus obtained. Similarly, we retain the notation for various other objects associated with our family of curves.

\subsection
\label{HnNotation}
{\it Notation.\/} --- Let $M$ be a module over some commutative $R$-algebra, or a sheaf on some $R$-scheme, on which the group scheme~$\mu_m$ acts. For $n \in \mZ/m\mZ$ we write
$$
M_{(n)} := \bigl\{x \in M \bigm| \zeta(x) = \zeta^n \cdot x \bigr\}\, ,
$$
which in the sheaf case has to be interpreted on the level of local sections. We refer to~$M_{(n)}$ as the $n$-eigenspace of~$M$. We have $M = \oplus_{n \in \mZ/m\mZ}\, M_{(n)}$.

\subsection
Recall that $r_i = \gcd(m,a_i)$. Consider the $\mu_m$-cover $\pi \colon C_t \to \mP^1$ for some $t \in T(k)$, where $k$ is a field. For $i\in \{1,\ldots,N\}$ and $n \in \mZ$, let 
$$
l(i,n) := -1 + \left\lceil {r_i - n a_i \over m}\right\rceil  = \left\lfloor {-n a_i\over m}\right\rfloor\, ,
$$
where the second equality easily follows from the fact that $r_i = \gcd(m,a_i)$. Consider the differential forms
$$
\omega_{n,\nu} := y^n \cdot (x-t_1)^\nu \cdot \prod_{i=1}^N (x-t_i)^{l(i,n)} \cdot dx\, ,\eqlabel{omegannu}
$$
and note that these only depend on the pair $(n \bmod m,\nu)$.

The following result is standard.

\subsection
\label{DiffFormsLemma}
{\it Lemma. --- Let $n \in \mZ/m\mZ$ with $n \neq 0$. The forms $\omega_{n,\nu}$ for $0 \leq \nu \leq  -2 + \sum_{i=1}^N\, \frpart{{-n a_i\over m}}$ are regular $1$-forms on~$C_t$ and they form a $k$-basis for $H^0(C_t,\Omega^1)_{(n)}$.\/}

\section{The special subvariety given by the action of the covering group}{SpecialSubvar}
\bigskip

\noindent
In this section we work over~$\mC$ and we fix a triple $(m,N,a)$ as in~\ref{CoverTsetup} with $N \geq 4$. We retain the notation introduced in the previous section, with the convention that all objects are now considered over~$\mC$ via the chosen embedding $R \hookrightarrow \mC$.

\subsection
{\it Definition.\/} --- We let $Z(m,N,a) \subset \cA_{g,\mC}$ be the scheme-theoretic (or rather, stack-theoretic) image of the morphism $\phi \colon T \to \cA_g$ of~(\ref{phiDef}). 
\medskip

In other words, $Z(m,N,a)$ is the reduced closed substack of~$\cA_{g,\mC}$ with underlying topological space $\ol{\phi(T)}$. We note that $Z(m,N,a)$ only depends on the equivalence class of $(m,N,a)$ and does not depend on the choice of the open subscheme~$T$ in~\ref{CoverTsetup}. The dimension of $Z(m,N,a)$ equals $N-3$.

\subsection
{\it Notation.\/} --- It will be convenient to write
$$
I(m) := \bigl[(\mZ/m\mZ)\setminus \{0 \bmod m\}\bigr]/\{\pm 1\}\, .
\eqlabel{I(m)def}
$$
Recall that $\frpart{x}$ denotes the fractional part of~$x$. For $n \in \mZ/m\mZ$ we define
$$
d_n := \cases{-1 + \sum_{i=1}^N\, \frpart{{-n a_i\over m}} & if $n \not\equiv 0$,\cr 0 & if $n \equiv 0$,}
\eqlabel{dndef}
$$
which by Lemma~\ref{DiffFormsLemma} is the dimension of the $(n)$-eigenspace of $H^0(C_t,\Omega^1)$, for any $t \in T(\mC)$. 

\subsection
\label{SmumDef}
The substack $Z(m,N,a) \subset \cA_g$ is contained in a special subvariety $S(\mu_m) \subset \cA_g$ determined by the action of~$\mZ[\mu_m]$ on the relative Jacobian $J \to T$. More precisely, $S(\mu_m)$ is the largest closed, reduced and irreducible substack $S\subset \cA_g$ containing~$Z(m,N,a)$ such that the homomorphism $\mZ[\mu_m] \to \End(J/T)$ induced by the action of~$\mu_m$ on $C/T$ extends to an action of $\mZ[\mu_m]$ on the universal abelian scheme over~$S$.

Choose a base point $b \in T(\mC)$, and let $(J_b,\lambda)$ be the corresponding Jacobian with its principal polarization. With $(V_\mZ,\Psi)$ as in~\ref{SiegelMod}, choose a symplectic similitude $\sigma\colon H_1(J_b,\mZ) \isomarrow V_\mZ$, where we equip $H_1(J_b,\mZ)$ with its Riemann form (i.e., the polarization, in the sense of Hodge theory, that corresponds with~$\lambda$). Via~$\sigma$, the action of~$\mu_m$ on~$J_b$ induces a structure of a $\mQ[\mu_m]$-module on $V = V_\mZ \otimes \mQ$. Consider the algebraic subgroup $H \subset G_\mQ = \CSp(V,\Psi)$ of $\mQ[\mu_m]$-linear symplectic similitudes, i.e., the subgroup given by
$$
H := \GL_{\mQ[\mu_m]}(V) \cap \CSp(V,\Psi)\, .
$$
With notation as in~\ref{SiegelMod}, the image of $Y_H \subset \gH$ under the map
$$
\gH \twoheadrightarrow G(\mZ)\backslash \gH \cong G(\mQ)\backslash \gH \times G(\Af)/G(\hat{\mZ}) \cong \cA_g(\mC)
$$
is (the set of $\mC$-points of) a finite union of algebraic subvarieties of~$\cA_g$, and the special subvariety~$S(\mu_m)$ is the unique irreducible component of this image that contains~$Z(m,N,a)$.

The dimension of~$S(\mu_m)$ is given by
$$
\dim S(\mu_m) = \sum\, d_{-n}d_n + \cases{{d_k(d_k+1)\over 2} & if $m=2k$ is even,\cr 0 & if $m$ is odd,\cr}
\eqlabel{dimSmum}
$$
where the first sum runs over the pairs $\pm n \in I(m)$ with $2n \not\equiv 0$. See [\ref{BMFO}], Section~5; see also Remark~\ref{SmumDim} below.

In particular, as $Z(m,N,a) \subset S(\mu_m)$ we have
$$
N-3 \leq \dim S(\mu_m)\, .\eqlabel{diminequality}
$$
If equality holds then $Z(m,N,a) = S(\mu_m) \subset \cA_g$ is a special subvariety; it then follows that among the Jacobians~$J_t$, for $t \in T(\mC)$, there are, up to isomorphism, infinitely many Jacobians of CM type. (Here we use that on a special subvariety the CM points lie dense.)

\subsection
{\it Inventory of examples.\/} --- Table~\ref{tablespecials} lists twenty triples $(m,N,a)$ with $N \geq 4$ for which $N-3 = \dim S(\mu_m)$, so that $Z(m,N,a) \subset \cA_{g,\mC}$ is a special subvariety. Our assumption that $N\geq 4$ means we are only considering the cases that give rise to a special subvariety of positive dimension.

The first column of the table gives a number that we assign to each example for reference. The examples are sorted first by genus, then by degree of the cover. The second column gives the genus of the curves in the family. In the next three columns we give the data $(m,N,a)$. In most cases there is a unique $a = (a_1,\ldots,a_N)$ in its equivalence class such that $1 \leq a_1 \leq \cdots \leq a_N \leq m-1$ with $\sum a_i = 2m$, and if there is a unique such representative, this is the one we list. If $N=4$ there are usually two such representatives, the second one being $(m-a_4,m-a_3,m-a_2,m-a_1)$; we list the one which lexicographically comes first. In the last column we give references to places in the literature where the example can be found. Here ``S($i$)'' refers to example~($i$) in~[\ref{Shimura}], ``M$i$'' refers to example~$i$ in the appendix of~[\ref{Most}], and dJN(1.3.$i$) refers to example~(1.3.$i$) in~[\ref{dJN}]. These examples can also be found in~[\ref{Rohde}]; see \ref{RohdeThm} below.

\topinsert
\tablelabel{tablespecials}
$$
\vbox{
\def\nr{\global\advance\caseno by 1({\oldstyle\casenr})}
\def\horrule{\noalign{\smallskip}\noalign{\hrule}\noalign{\medskip}}
\halign{\hfill #\hfill\quad& \quad\hfil$#$\hfil\quad & \quad\hfil$#$\hfil\quad & \quad\hfil$#$\hfil\quad & \quad\hfil$#$\hfil\quad & \qquad$#$\hfil\quad \cr
\omit \hfill & \hbox{genus} & m & N & a & \hbox{references} \cr
\horrule
\nr & 1 & 2 & 4 & (1,1,1,1) & \cr
\horrule
\nr & 2 & 2 & 6 & (1,1,1,1,1,1) & \cr
\horrule
\nr & 2 & 3 & 4 & (1,1,2,2) & \cr
\horrule
\nr & 2 & 4 & 4 & (1,2,2,3) & \cr
\horrule
\nr & 2 & 6 & 4 & (2,3,3,4) & \cr
\horrule
\nr & 3 & 3 & 5 & (1,1,1,1,2) & \hbox{S(1), M41} \cr
\horrule
\nr & 3 & 4 & 4 & (1,1,1,1) & \cr
\horrule
\nr & 3 & 4 & 5 & (1,1,2,2,2) & \hbox{S(3), M43} \cr
\horrule
\nr & 3 & 6 & 4 & (1,3,4,4) & \cr
\horrule
\nr & 4 & 3 & 6 & (1,1,1,1,1,1) & \hbox{S(2), M23, dJN(1.3.1)} \cr
\horrule
\nr & 4 & 5 & 4 & (1,3,3,3) & \hbox{S(4), dJN(1.3.2)} \cr
\horrule
\nr & 4 & 6 & 4 & (1,1,1,3) & \cr
\horrule
\nr & 4 & 6 & 4 & (1,1,2,2) & \cr
\horrule
\nr & 4 & 6 & 5 & (2,2,2,3,3) & \cr
\horrule
\nr & 5 & 8 & 4 & (2,4,5,5) &\cr
\horrule
\nr & 6 & 5 & 5 & (2,2,2,2,2) & \hbox{S(5), M44} \cr
\horrule
\nr & 6 & 7 & 4 & (2,4,4,4) & \hbox{S(6), dJN(1.3.3)} \cr
\horrule
\nr & 6 & 10 & 4 & (3,5,6,6) & \cr
\horrule
\nr & 7 & 9 & 4 & (3,5,5,5) & \cr
\horrule
\nr & 7 & 12 & 4 & (4,6,7,7) & \cr
\horrule\cr}}
$$
\centerline{{\bf Table~\ref{tablespecials}: monodromy data that give rise to a special subvariety}}
\bigskip
\noindent\hrule
\endinsert

\subsection
{\it Remark.\/} --- If we have a triple $(m,N,a)$ as in~\ref{CoverTsetup} with $N=3$, the associated subvariety $Z= Z(m,N,a)$ is a special point in~$\cA_g$.
\medskip

The following theorem is the main result of this paper. It says that the list of examples in Table~\ref{tablespecials} is exhaustive.

\subsection
\label{mainthm}
{\it Theorem. --- Consider data $(m,N,a)$ as in\/~{\rm \ref{CoverTsetup}}, with $N \geq 4$. Then $Z(m,N,a) \subset \cA_{g,\mC}$ is a special subvariety if and only if $(m,N,a)$ is equivalent to one of the twenty examples listed in\/~{\rm Table~\ref{tablespecials}}.}
\medskip

By what was explained above, it only remains to be shown that $Z(m,N,a)$ is not special if $(m,N,a)$ is not equivalent to one of the triples in the table. In~\ref{pfm=2} we shall first prove this for $m=2$. In Section~\ref{PfN=4} we shall prove the theorem for $N=4$. The case $N>4$ is treated in Section~\ref{PfN>4}.

Combining the theorem with the main result of~[\ref{Yafaev}] we obtain, for $N=4$, the following finiteness result for the number of CM fibres.

\subsection
\label{CorMain}
{\it Corollary. --- Assume the Generalized Riemann Hypothesis for CM fields. If we have a triple $(m,N,a)$ with $N=4$ and the equivalence class of this triple is not among the twenty examples listed in\/~{\rm Table~\ref{tablespecials}} then up to isomorphism there are finitely many Jacobians of CM type in the family $J \to T_\mC$.}

\section{Decomposition of the variation of Hodge structure}{DecompVHS}
\bigskip

\noindent
In this section we again work over~$\mC$. We fix a triple $(m,N,a)$ as in~\ref{CoverTsetup} with $N \geq 4$. The corresponding family of curves $C \to T$ gives rise to a variation of Hodge structure over~$T$ (over~$\mC$) on which $\mQ[\mu_m]$ acts. We describe the resulting decomposition of the VHS, and we give some results about the monodromy of the summands.

\subsection
\label{reduction}
Let $(m,N,a)$ be as in~\ref{CoverTsetup}. For $n \in \mZ/m\mZ$, let $N(n)$ be the number of indices $i \in \{1,\ldots,N\}$ such that $n a_i \not\equiv 0$ modulo~$m$. In particular, $N(n) = N$ if $n \in (\mZ/m\mZ)^*$. If $n \neq 0$ then $d_n + d_{-n} = N(n) - 2$.

Let $m^\prime \geq 2$ be a divisor of~$m$, say with $m = rm^\prime$. Write $N^\prime := N(r \bmod m)$, and consider the $N^\prime$-tuple $a^\prime$ in $(\mZ/m^\prime\mZ)^{N^\prime}$ obtained from the $N$-tuple $(a_1 \bmod {m^\prime},\ldots,a_N \bmod {m^\prime})$ by omitting the zero entries. Then $(m^\prime,N^\prime,a^\prime)$ is again a triple as in~\ref{CoverTsetup}; we refer to it as the triple obtained from $(m,N,a)$ by reduction modulo~$m^\prime$.

\subsection
\label{monodromy}
As in~\ref{FamilyConstr}, consider the family of curves $f \colon C \to T$ associated with the triple $(m,N,a)$, where, as in the previous section, we work over~$\mC$. The cohomology groups $H^1(C_t,\mQ)$ are the fibres of a polarized $\mQ$-VHS with underlying local system $\mV := R^1 f^\an_* \mQ_C$. This $\mQ$-VHS comes equipped with an action of the group ring $\mQ[\mu_m] = \prod_{d|m} K_d$, where $K_d = \mQ[t]/\Phi_d$ is the cyclotomic field of $d$th roots of unity. Accordingly we have a decomposition of $\mQ$-VHS,
$$
\mV = \oplus_{d|m} \mV^{[d]}\, ,
\eqlabel{V[d]Dec}
$$ 
where $\mV^{[d]}$ is a polarized $\mQ$-VHS over~$T$ equipped with an action of~$K_d$.

\subsection
\label{oldfactors}
Let $(m,N,a)$ be as in~\ref{CoverTsetup}, with $N \geq 4$. Let $m^\prime \geq 2$ be a proper divisor of~$m$, and let $(m^\prime,N^\prime,a^\prime)$ be the triple obtained from $(m,N,a)$ by reduction modulo~$m^\prime$, as in~\ref{reduction}. 

Let $\rho \colon \mA^N \to \mA^{N^\prime}$ be the projection map, omitting the coordinates~$x_i$ for all indices~$i$ with $a_i \equiv 0$ modulo~$m^\prime$. Performing the construction of~\ref{FamilyConstr} we may choose $T = T(m,N,a) \subset \mA^N$ and $T^\prime = T(m^\prime,N^\prime,a^\prime) \subset \mA^{N^\prime}$ such that $\rho$ maps~$T$ to~$T^\prime$. (The map $T \to T^\prime$ is then dominant.) If $f^\prime \colon C^\prime \to T^\prime$ is the family of curves associated with the triple $(m^\prime,N^\prime,a^\prime)$, the pull-back $\rho^* C^\prime = C^\prime \times_{T^\prime} T$ can be identified with the quotient of~$C/T$ modulo the action of $\mu_r \subset \mu_m$, where $r = m/m^\prime$. 

Let $\mV^\prime = \mV(m^\prime,N^\prime,a^\prime)$ be the $\mQ$-VHS over~$T^\prime$ associated with the triple $(m^\prime,N^\prime,a^\prime)$. The quotient morphism $C \to \rho^* C^\prime$ over~$T$ gives a map $\rho^* \mV^\prime \to \mV$ and this induces an isomorphism $\rho^* \mV^\prime \cong \oplus_{d|m^\prime} \mV^{[d]}$. We refer to the sub-$\mQ$-VHS $\rho^* \mV^\prime$ as the {\it old part\/} associated with the divisor~$m^\prime$. The direct summand $\mV^{[m]} \subset \mV$ is called the {\it new part\/}; in the decomposition~(\ref{V[d]Dec}) it is the only direct summand that is not contained in an old part.

\subsection
\label{oldspecial}
{\it Lemma. --- Let the notation and assumptions be as in\/~{\rm \ref{oldfactors}}. Let $g = g(m,N,a)$ and $g^\prime = g(m^\prime,N^\prime,a^\prime)$ be the respective genera. If $Z(m,N,a) \subset \cA_{g,\mC}$ is a special subvariety then  $Z(m^\prime,N^\prime,a^\prime) \subset \cA_{g^\prime,\mC}$ is special, too.}
\medskip

\Proof
Suppose $Z(m,N,a)$ is special. Write $J = J(m,N,a) \to T$ and $J^\prime = J(m^\prime,N^\prime,a^\prime) \to T^\prime$ for the respective Jacobians. Then $\rho^* J^\prime$ is an isogeny factor of~$J$. Hence, if $\theta\colon T \to \cA_{g^\prime,\mC}$ is the morphism corresponding to~$\rho^* J^\prime$, it follows from~\ref{IsogenyFactors} that the Zariski closure of the image of~$\theta$ is a special subvariety. But $\theta$ is just the composition of the projection $\rho \colon T \to T^\prime$, which is dominant, and the morphism $\phi(m^\prime,N^\prime,a^\prime) \colon T^\prime \to \cA_{g^\prime,\mC}$ associated with the data $(m^\prime,N^\prime,a^\prime)$. Hence the closure of the image of~$\theta$ is $Z(m^\prime,N^\prime,a^\prime)$. 
\QED

\subsection
With notation as in~\ref{HnNotation}, the $\mC$-local system~$\mV_\mC$ decomposes as $\mV_\mC = \oplus_{n \in \mZ/m\mZ} \mV_{\mC,(n)}$. The relation with~(\ref{V[d]Dec}) is that $\mV^{[d]}_\mC$ is the sum of all $\mV_{\mC,(n)}$ with $\order(n) = d$. We have $V_{\mC,(0)} = 0$, and for $n \neq 0$ the $\mC$-local system $\mV_{\mC,(n)}$ has rank $N(n)-2$, where $N(n)$ is defined as in~\ref{reduction}. (See also [\ref{DelMost}], Section~2.)

With real coefficients, and with $I(m)$ as in~(\ref{I(m)def}), we have a decomposition of $\mR$-VHS
$$
\mV_\mR = \left(\bigoplus_{{\pm n \in I(m), 2n \not\equiv 0}}\, \mV_{\mR,(\pm n)} \right) \oplus \mV_{\mR,({m\over 2})}\, ,
$$
where the last summand only occurs if $m$ is even. The decomposition is such that
$$
\mV_{\mR,(\pm n)} \otimes_\mR \mC = \mV_{\mC,(n)} \oplus \mV_{\mC,(-n)}
$$
and
$$
\mV_{\mR,({m\over 2})} \otimes_\mR \mC = \mV_{\mC,({m\over 2})}\quad \hbox{(for even $m$)}\, .
$$
On the summand $\mV_{\mR,(\pm n)}$ the polarization of the VHS induces a $(-1)$-hermitian form~$\beta_{\pm n}$ of signature $(d_n,d_{-n})$; see also [\ref{DelMost}], Corollary~(2.21) and~(2.23). In case $m$ is even, the polarization induces a symplectic form $\beta_{{m\over 2}}$ on $\mV_{\mR,({m\over 2})}$.

As in \ref{SmumDef} we choose a base point $b \in T(\mC)$, and a symplectic similitude $\sigma \colon H^1(J_b,\mQ) \isomarrow V$. Via this similitude, we identify~$V$ with the fibre of~$\mV$ at~$b$. Correspondingly, we write $V_{\mR,(\pm n)}$ for the direct summand of~$V_\mR$ that under~$\sigma$ maps to the fibre at~$b$ of $\mV_{\mR,(\pm n)}$. 

For what follows, it will be convenient to choose the base point~$b$ to be a Hodge-generic point with respect to the variation~$\mV$, so from now on we assume this. Define $\Mono \subset \GL(V)$ and $\Hdg \subset \GL(V)$ to be the algebraic monodromy groups of the $\mQ$-local system~$\mV$ and the Hodge group of the fibre at~$b$, respectively. (Since we have chosen~$b$ to be Hodge-generic, $\Hdg$ is the generic Hodge group of the VHS.) By [\ref{Andre}], Theorem~1, the identity component $\Mono^0$ is a semisimple normal subgroup of~$\Hdg$.

\subsection
\label{SmumDim}
{\it Remark.\/} --- The generic Hodge group~$\Hdg$ is contained in the algebraic group $H := \GL_{\mQ[\mu_m]}(V) \cap \CSp(V,\Psi)$ considered in~\ref{SmumDef}. Extending scalars to~$\mR$ we have 
$$
H_\mR \cong \left(\prod_{{\pm n \in I(m), 2n \not\equiv 0}}\, \UU(V_{\mR,(\pm n)},\beta_{\pm n}) \right) \times \Sp(V_{\mR,({m\over 2})},\beta_{{m\over 2}})\, ,
\eqlabel{HRDec}
$$
where the symplectic factor only occurs if $m$ is even. Our formula (\ref{dimSmum}) follows from this together with what was explained in~\ref{MTfactors}.
\medskip

We have $\Mono^0 \subset \Hdg \subset H$. The next proposition gives us information about the projections of the $\mR$-group~$\Mono^0_\mR$ onto the various factors in the decomposition~(\ref{HRDec}).

\subsection
\label{monosurjects}
{\it Proposition. --- {\rm (\romno1)} Suppose we have $\pm n \in I(m)$ with $2n \not\equiv 0$ such that $d_n \geq 1$ and $d_{-n} \geq 1$. Then the image of $\Mono^0_\mR$ in $\UU(V_{\mR,(\pm n)},\beta_{\pm n})$ is the special unitary group $\SU(V_{\mR,(\pm n)},\beta_{\pm n})$, which is isomorphic to $\SU(d_n,d_{-n})$.

{\rm (\romno2)} If $m$ is even, the group $\Mono^0_\mR$ projects surjectively to the symplectic group $\Sp(V_{\mR,({m\over 2})},\beta_{{m\over 2}})$, which is isomorphic to $\Sp_{2h,\mR}$ with $h = d_{{m\over 2}}$.}
\medskip

\Proof
For (\romno1), see Rohde~[\ref{Rohde}], Theorem~5.1.1.  For~(\romno2), assume $m$ is even, and let $h := d_{{m\over 2}}$. There are then $2h+2$ indices~$i$ for which~$a_i$ is even, and without loss of generality we may assume these are the indices $i \in \{1,\ldots,2h+2\}$. With the above notation, $\mV_{\mR,({m\over 2})}$ is the same as $\mV^{[2]}_\mR$, and if we define $\Mono^{[2]} \subset \GL(V^{[2]})$ as the algebraic monodromy group of the $\mQ$-local system~$\mV^{[2]}$ then the image of~$\Mono^0$ in $\Sp(V_{\mR,({m\over 2})},\beta_{{m\over 2}})$ equals $\Mono^{[2],0}_\mR$. Further, $\mV^{[2]}$ is just the local system obtained from the family of hyperelliptic curves $u^2 = (x-t_1)\cdots(x-t_{2h+2})$; cf.~\ref{oldfactors}. By [\ref{ACampo}], Theorem~1 (see also [\ref{AchtPries}]), the connected algebraic monodromy group of this local system is the full symplectic group. 
\QED

\subsection
\label{RohdeThm}
Let $(m,N,a)$ be a triple as in~\ref{CoverTsetup}, with $N \geq 4$. Consider the following condition.
$$
\vcenter{
\vskip 2pt
\setbox0=\hbox{There exists a $\pm n \in I(m)$ with $\{d_n,d_{-n}\} = \{1,N-3\}$,}
\copy0
\vskip 1pt
\hbox to\wd0{\hfill and for all other $\pm n^\prime \in I(m)$, either $d_{n^\prime}=0$ or $d_{-n^\prime} =0$.\hfill}}
\eqlabel{RohdeCond}
$$

All examples in Table~\ref{tablespecials}, with the exception of Example~{\oldstyle 2}, satisfy this condition. It has been proven by Rohde in~[\ref{Rohde}] that these nineteen examples are the only examples (up to equivalence) of triples $(m,N,a)$ with $N \geq 4$ that satisfy~(\ref{RohdeCond}).

\subsection
\label{pfm=2}
We now prove Theorem~\ref{mainthm} for $m=2$. In this case $N$ is even, $g=(N-2)/2$, and up to equivalence $a=(1,\ldots,1)$ is the only possibility. Further, $Z(2,N,a) \subset \cA_{g,\mC}$ is the closure of the hyperelliptic locus in~$\cA_g$. The assertion is that this is not special for $g\geq 3$. This is true because for $g\geq 3$ the hyperelliptic locus is not dense in~$\cA_g$, whereas by (\romno2) of Proposition~\ref{monosurjects} the generic Hodge group~$\Hdg$ is the full symplectic group.
\medskip

In what follows we may therefore assume $m>2$.

\section{Canonical liftings of ordinary Jacobians in positive characteristic}{canlifts}

\subsection
\label{KS}
We need to recall some definitions and results from the paper [\ref{DworkOgus}] by Dwork and Ogus. We first discuss their theory in a general setting. After that, from \ref{assumpt} on, we return to the specific families of curves that interest us, and we explain how the results of~[\ref{DworkOgus}] apply to that situation.

Given an irreducible base scheme~$T$ and a smooth projective curve $f \colon C \to T$ we write $\omega_{C/T}$ for $\Omega^1_{C/T}$. Consider the Hodge bundle $\mE = \mE(C/T) := f_* \omega_{C/T}$, which is a vector bundle on~$T$ of rank~$g$, the genus of the fibres. We have a Kodaira-Spencer map $\KS\colon \Sym^2(\mE) \to \Omega^1_{T/R}$ and we define
$$
\eqalign{
&\cK = \cK(C/T) := \Ker\bigl(\Sym^2(\mE) \mapright{\mult} f_*(\omega_{C/T}^{\otimes 2})\bigr)\, ,\cr
&\cQ = \cQ(C/T) := \Coker\bigl(f_*(\omega_{C/T}^{\otimes 2})^\vee \mapright{\mult^\vee} \Sym^2(\mE)^\vee\bigr)\, ,\cr}
$$
where ``$\mult$'' is the multiplication map. We remark that $\cQ$ can be thought of as the pull-back to~$T$ of the normal bundle to the Torelli locus inside~$\cA_g$. Outside the hyperelliptic locus, $\mult$ is surjective and $\cK$ is the dual of~$\cQ$. 

\subsection
\label{DO}
Let $k$ be an algebraically closed field of characteristic $p>0$. Consider a smooth projective curve $C/k$ of genus~$g$ such that its Jacobian~$J$ is ordinary. Let $\lambda$ be the natural principal polarization of~$J$. By Serre-Tate theory (see~[\ref{Katz}] or [\ref{Messing}], Chap.~5) we have a canonical lifting of $(J,\lambda)$ to a principally polarized abelian variety $(J^\can,\lambda^\can)$ over the ring of Witt-vectors~$W(k)$. In general, for $g \geq 4$ this canonical lifting is not the Jacobian of a curve over~$W(k)$. More precisely, Dwork and Ogus show that in general not even the canonical lifting over~$W_2(k)$, the Witt-vectors of length~$2$, is again a Jacobian.

Following~[\ref{DworkOgus}] we call the ordinary curve~$C/k$ {\it pre-$W_2$-canonical\/} if there exists a smooth projective curve~$Y$ over~$W_2(k)$ such that $(J^\can,\lambda^\can)$ over~$W_2(k)$ is isomorphic, as a principally polarized abelian variety, to the Jacobian of~$Y$.

The main advantage of working modulo~$p^2$ is that in this case there is a natural class $\beta_{C/k}$ in $\cQ(C/k)$ that vanishes if and only if $C/k$ is pre-$W_2$-canonical; see [\ref{DworkOgus}], Prop.~(2.4). This invariant~$\beta_{C/k}$ cannot be defined in a family of curves (in a sense that can be made precise, see [\ref{DworkOgus}], \S~3), but its Frobenius pullback can, as we shall discuss next.

Let $T$ be a smooth $k$-scheme, and let $\Frob_T \colon T \to T$ be its absolute Frobenius endomorphism. Consider a smooth projective curve $f \colon C \to T$ such that all fibres~$C_t$ are ordinary. This assumption permits us to view the inverse Cartier operator as an $\cO_T$-linear homomorphism 
$$
\gamma \colon \Frob_T^* \mE \to \mE\, ,
$$
with $\mE = \mE(C/T)$ the Hodge bundle. This map~$\gamma$ is the inverse transpose of the Frobenius action on $R^1 f_* \cO_C$. 

Dwork and Ogus define a global section $\tilde\beta_{C/T}$ of $\Frob_T^* \cQ(C/T)$ such that for any $t \in T(k)$ the value of~$\tilde\beta_{C/T}$ at~$t$ is $\Frob_k^*(\beta_{C_t/k})$. See [\ref{DworkOgus}], Prop.~2.7. 

The sheaf $\Frob_T^* \cQ(C/T)$ comes equipped with a canonical flat connection
$$
\nabla \colon \Frob_T^* \cQ(C/T) \to \Frob_T^* \cQ(C/T) \otimes \Omega^1_{T/k}\, .
$$
In particular, still with $T$ smooth over~$k$, if $C_t/k$ is pre-$W_2$-canonical for all $t \in T(k)$, the section~$\tilde\beta_{C/T}$ is zero, so certainly $\nabla \tilde\beta_{C/T} = 0$. One of the main points of~[\ref{DworkOgus}], then, is that $\nabla \tilde\beta_{C/T}$ can be calculated explicitly. The result is easiest to state under the additional hypothesis that the fibres $C_t$ are not hyperelliptic. We note that this assumption is not always satisfied in the situation we want to consider. In Prop.~\ref{Criterion} we shall give a modified version of the next result that also works for hyperelliptic families. 

If the $C_t$ are not hyperelliptic, the multiplication map $\Sym^2(\mE) \to f_*(\omega^{\otimes 2})$ is surjective, so, with notation as introduced in~\ref{KS}, $\tilde\beta_{C/T}$ can be viewed as an $\cO_T$-linear map $\Frob_T^* \cK \to \cO_T$.

The following result is [\ref{DworkOgus}], Prop.~3.2. (In~[\ref{DworkOgus}] the result is stated under some versality assumption on the family~$C/T$ but one easily verifies that their Prop.~3.2 is true without this assumption.)

\subsection
\label{nablabeta}
{\it Proposition (Dwork-Ogus). --- With assumptions and notation as above, and assuming the curves~$C_t$ to be non-hyperelliptic, $-\nabla(\tilde\beta_{C/T}) \colon \Frob_T^* \cK \to \Omega^1_{T/k}$ is equal to the composition 
$$
\Frob_T^* \cK \injmapright{{\rm incl}} \Frob_T^* \Sym^2(\mE) \mapright{\Sym^2(\gamma)} \Sym^2(\mE) \mapright{\KS} \Omega^1_{T/k}\, . 
$$}
\vskip-\lastskip

\subsection
\label{assumpt}
{\it Assumptions.\/} --- In the rest of this section we fix data $(m,N,a)$ as in~\ref{CoverTsetup} with $m>2$ and $N \geq 4$, and we consider the family of curves $f \colon C \to T$ as constructed in~\ref{FamilyConstr}. We assume that $(m,N,a)$ is not equivalent to one of the triples listed in Table~1. We further assume that the closed subvariety $Z(m,N,a) \subset \cA_{g,\mC}$ is a special subvariety; our aim is to derive a contradiction. 

Without loss of generality we may suppose $(m,N,a)$ is minimal with respect to the above assumptions, by which we mean that there is no triple $(m^\prime,N^\prime,a^\prime)$ that satisfies our assumptions and for which we have $4 \leq N^\prime < N$ or $N^\prime = N$ and $2 \leq m^\prime < m$. In particular, it follows from Lemma~\ref{oldspecial} that if $m^\prime$ is a proper divisor of~$m$ and $(m^\prime,N^\prime,a^\prime)$ is obtained from $(m,N,a)$ by reduction modulo~$m^\prime$, as in~\ref{oldfactors}, either $N^\prime \leq 3$ or $(m^\prime,N^\prime,a^\prime)$ is equivalent to one of the twenty triples listed in Table~\ref{tablespecials}.

\subsection
\label{KSlemma}
{\it Lemma. --- {\rm (\romno1)} The subsheaf $f_*(\omega_{C/T}^{\otimes 2})_{(0)}$ of $\mu_m$-invariant sections of $f_*(\omega_{C/T}^{\otimes 2})$ is a locally free $\cO_T$-module of rank~$N-3$.

{\rm (\romno2)} The Kodaira-Spencer morphism $\KS \colon \Sym^2(\mE) \to \Omega^1_{T/R}$ factors over the composite map 
$$
\Sym^2(\mE) \mapright{\pr} \Sym^2(\mE)_{(0)} \mapright{\mult} f_*(\omega_{C/T}^{\otimes 2})_{(0)}
$$ 
and the induced map $\KS_{(0)} \colon f_*(\omega_{C/T}^{\otimes 2})_{(0)} \to \Omega^1_{T/R}$ is a fibrewise injective homomorphism of vector bundles.}
\medskip

\Proof
(\romno1) By a result of Chevalley and Weil in~[\ref{ChevWeil}], for every $t \in T(\mC)$ the subspace of $\mu_m$-invariants in $H^0(C_t,\Omega^{\otimes 2})$ has dimension $N-3$. 

(\romno2) The fibre of $\Sym^2(\mE)_{(0)}$ at a point~$t$ is dual to the space $H^1(C_t,\Theta)_{(0)}$, which parametrizes the deformations of~$C_t$ for which the $\mu_m$-action deforms along. Since we have a $\mu_m$-action on the whole family $C/T$, the Kodaira-Spencer map factors through the projection $\Sym^2(\mE) \to \Sym^2(\mE)_{(0)}$. For any family of curves it also factors through the multiplication map $\Sym^2(\mE) \to f_*(\omega_{C/S}^{\otimes 2})$; hence $\KS$ factors through $f_*(\omega_{C/T}^{\otimes 2})_{(0)}$. The last assertion of~(\romno2) just says that at any $t \in T(k)$ our family is complete in the sense of deformation theory: if $D \to \Spec\bigl(k[\epsilon]\bigr)$ is a first-order deformation of~$C_t$ with its $\mu_m$-action, $D/k[\epsilon]$ can be obtained by pull-back from our family. This is clear, as the quotient $D/\mu_m$ is isomorphic to~$\mP^1_{k[\epsilon]}$. 
\QED

\subsection
\label{preW2lemma}
{\it Lemma. --- There exists a prime number $p$ with $p \equiv 1$ modulo~$m$ and, for $\gp$ a prime of $R = \mZ[1/m,u]/\Phi_m$ above~$p$, a dense open subset $U$ of $T_0 := T \otimes_R (R/\gp)$, such that for any algebraically closed field~$k$ of characteristic~$p$ and any $t \in U(k)$ the curve $C_t$ is ordinary and the Jacobian~$J_t$ is pre-$W_2$-canonical.}
\medskip

\Proof
Let $p$ be any prime number with $p \equiv 1 \bmod m$. Note that in this case $R/\gp = \mF_p$ for any prime $\gp \subset R$ above~$p$. We write $T \otimes \mF_p$ for $T \otimes_R (R/\gp)$. By [\ref{Bouw}], Prop.~7.4, there is a Zariski open $U \subset T \otimes \mF_p$ such that $C_t$ is ordinary for all $t \in U(k)$. On the other hand, by a theorem of Noot~[\ref{NootThesis}], [\ref{NootModels}], in the slightly more precise formulation of~[\ref{BMLin2}], Thm.~4.2, the assumption that $Z= Z(m,N,a)$ is special implies that for $p$ large enough and $t \in T_0(k)$ any ordinary point, the canonical lifting of~$J_t$ gives a $W(k)$-valued point of~$Z$. In particular, $J_t$ is then pre-$W_2$-canonical.
\QED

\subsection
\label{pUChoice}
For most choices of $(m,N,a)$ the general member in our family of curves $C \to T$ is not hyperelliptic. In this case we may, and will, choose $p$ and~$U$ in Lemma~\ref{preW2lemma} such that the curves $C_t$ with $t \in U(k)$ are non-hyperelliptic.

It may happen, however, that $Z(m,N,a)$ is fully contained in the hyperelliptic locus, in which case we call $(m,N,a)$ a hyperelliptic triple. This occurs, for instance, if $m=2m^\prime$ is even and $a$ is of the form $(a_1,a_2,m^\prime,\ldots,m^\prime)$, or if $N=4$ and $a = (1,1,m-1,m-1)$. In such a case we will just choose $p$ and~$U$ as in Lemma~\ref{preW2lemma}.

The following variant of Prop.~\ref{nablabeta} is the essential tool in our proof of Theorem~\ref{mainthm}. Let $p$ and~$U$ be chosen as above, and let $C_U \to U$ be the restriction of~$C/T$ to~$U$. With notation as in~\ref{HnNotation} and~\ref{KS} we have eigenspace decompositions $\cQ(C_U/U) = \oplus_{n \in \mZ/m\mZ}\, \cQ_{(n)}$; likewise for other $\cO_U[\mu_m]$-modules such as $\cK(C_U/U)$ and $\mE_U = \mE(C_U/U)$. Accordingly we can write $\tilde\beta_{C_U/U} = \sum_n\, \tilde\beta_{(n)}$, with $\tilde\beta_{(n)}$ a section of~$F_U^* \cQ_{(n)}$.

\subsection
\label{Criterion}
{\it Proposition. --- With assumptions and notation as in\/~{\rm \ref{assumpt}} and\/~{\rm \ref{pUChoice}}, the composite map
$$
\Frob_U^* \cK_{(0)} \injmapright{{\rm incl}} \Frob_U^* \Sym^2(\mE_U)_{(0)} \mapright{\Sym^2(\gamma)} \Sym^2(\mE_U)_{(0)} \mapright{\mult_{(0)}} f_*(\omega_{C/U}^{\otimes 2})_{(0)}
\eqlabel{composite(0)}
$$
is zero.}
\medskip

\Proof
If $(m,N,a)$ is non-hyperelliptic this is immediate from Prop.~\ref{nablabeta} together with (\romno2) of Lemma~\ref{KSlemma}. So from now on we may assume that $Z(m,N,a)$ is contained in the hyperelliptic locus. Let $\iota \in \Aut(C_U/U)$ be the hyperelliptic involution. We remark that $\iota$ may or may not be contained in the subgroup $\mu_m \subset \Aut(C_U/U)$; by inspection of the hyperelliptic examples mentioned in~\ref{pUChoice} we see that both cases occur. 

Let $\omega = \Omega^1_{C_U/U}$. In the hyperelliptic case the multiplication map $\Sym^2(\mE_U) \to f_*(\omega^{\otimes 2})$ is no longer surjective. However, if we denote the invariants under the action of~$\iota$ by a subscript ``even'',
$$
\mult_{(\even)} \colon \Sym^2(\mE_U)_{(\even)} \to f_*(\omega^{\otimes 2})_{(\even)}
$$
is again surjective; see for instance [\ref{FOJS}], Lemmas 2.12 and~2.13. The assumption that $Z(m,N,a)$ is fully contained in the hyperelliptic locus implies that $f_*(\omega^{\otimes 2})_{(0)} \subset f_*(\omega^{\otimes 2})_{(\even)}$. Hence also $\mult_{(0)}$ is surjective, and we may view $\tilde\beta_{(0)}$ as an $\cO_U$-linear map $F_U^* \cK_{(0)} \to \cO_U$. 

Our assumptions imply that $\tilde\beta_{(0)}$, hence also $\nabla\tilde\beta_{(0)}$, is zero, so all that remains to be checked is that for the $0$-component the analogue of Prop.~\ref{nablabeta} holds. For this we can follow the proof of [\ref{DworkOgus}], Prop.~3.2: The proof that loc.\ cit.\ diagram~(3.2.9) is commutative does not use the assumption that the curves are non-hyperelliptic and therefore goes through in the hyperelliptic case. Moreover, in the notation of that proof, we can choose the lifted family ${\bfit Y}/{\bfit T}$ (where $Y/T$ is our $C_U/U$) such that it has an action of~$\mu_m$. Taking $0$-components in diagram~(3.2.9) then gives that $-\nabla \tilde\beta_{(0)}$ equals the composition of~(\ref{composite(0)}) and~$\KS_{(0)}$, and by (\romno2) of Lemma~\ref{KSlemma} we conclude that the map~(\ref{composite(0)}) is zero.
\QED

\subsection
\label{InvgammaMatrix}
{\it Lemma. --- Write $q = (p-1)/m$. Let $n \in \mZ/m\mZ$, and let $A=A_n(t) \in \GL_{d_n}(\cO_U)$ be the matrix of the Cartier operator $\gamma^{-1}\colon \mE_{U,(n)} \to F_U^* \mE_{U,(n)}$ with regard to the frames 
$$
\omega_{n,0},\ldots,\omega_{n,d_n-1}\, ,\quad \hbox{respectively}\quad \omega_{n,0} \otimes 1,\ldots,\omega_{n,d_n-1} \otimes 1\, .
$$
Then the matrix coefficient $A_{\rho,\sigma}$, for $\rho$, $\sigma \in \{0,\ldots,d_n-1\}$, equals
$$
(-1)^\Sigma \cdot \sum_{j_1+\cdots+j_N=\Sigma}\, {q \cdot \repr{-na_1}_m \choose j_1} \cdots {q \cdot \repr{-na_N}_m \choose j_N} \cdot t_1^{j_1} \cdots t_N^{j_N}\, ,
\eqlabel{MatrCoeff}
$$
where $\Sigma = (d_n-\sigma)(p-1) + (\rho-\sigma)$.}
\medskip

\Proof
This is the dual version of [\ref{Bouw}], Lemma~5.1, part~(\romno1).
\QED
\bigskip

We remark that the result of Bouw in~[\ref{Bouw}] is valid in a more general context. She uses this to prove the existence of cyclic covers whose $p$-rank reaches a natural upper bound imposed by the local monodromy data; see [\ref{Bouw}], Thm.~6.1.

\section{Proof of the main result: four branch points}{PfN=4}
\bigskip

\noindent
In this section we prove Theorem~\ref{mainthm} in the case $N=4$. Let us first try to explain the idea of the argument, by way of guide for the reader. 

Assume we have a triple $(m,N,a)$ that is not one of the twenty examples in Table~\ref{tablespecials} but such that $Z(m,N,a)$ is special. We want to show that this is impossible. We have the family of curves $C \to T$ over a ring~$R$ of finite type over~$\mZ$; we then consider the reduction modulo~$p$, with~$p$ chosen as in~\ref{pUChoice}. The nontrivial information we have, exploiting the assumption that $Z(m,N,a)$ is special, is that this gives a family of ordinary Jacobians in characteristic~$p$ such that the canonical liftings over~$W_2(k)$ are again Jacobians. In Prop.~\ref{Criterion} we have translated this, based on the theory in~[\ref{DworkOgus}], into the vanishing of a certain abstractly defined map. Lemma~\ref{InvgammaMatrix} enables us to calculate this map explicitly. To apply this, we need a differential form that gives us an element in the source of the map~(\ref{composite(0)}); we write down such a form in~\ref{wpmncalc} below. The fact that the image of this element is zero then gives us a polynomial identity~(\ref{cc=cc}). By some combinatorial arguments we then show that this identity cannot hold if $(m,N,a)$ is not equivalent to one of the twenty special triples in Table~\ref{tablespecials}.

\subsection
We retain the assumptions made in~\ref{assumpt}. In addition we assume $N=4$. Let
$$
\cD(m,a) := \bigl\{\pm n \in I(m) \bigm| d_n = d_{-n} = 1\bigr\}\, .
$$
(Here $I(m)$ is as defined in~(\ref{I(m)def}) and $d_n$ is given by~(\ref{dndef}). Though not indicated in the notation, $d_n$ depends on~$a$.)

If $\pm n \in \cD(m,a)$, the equalities $d_n=d_{-n}=1$ imply that $na_i \not\equiv 0$ modulo~$m$ for all $i \in \{1,\ldots,4\}$. Hence,
$$
\repr{-na_i}_m = m - \repr{na_i}_m\quad \hbox{for all~$i$,}\, \qquad \hbox{and}\qquad \sum_{i=1}^4\, \repr{-n a_i}_m = 2m\, ,\eqlabel{-nai=m-nai}
$$
where we recall that $\repr{b}_m$ denotes the representative of the class $(b \bmod{m})$ in $\{0,1,\ldots,m-1\}$.

We first observe that $\dim S(\mu_m) > 1$. Indeed, if $\dim S(\mu_m) =1$ then (\ref{RohdeCond}) holds, and by the results of [\ref{Rohde}], Chapter~6, $(m,N,a)$ is one of the triples in Table~\ref{tablespecials}, which contradicts our assumptions.

For $n \in \mZ/m\mZ$ we have $d_n + d_{-n} \leq 2$, and if $m$ is even then $d_{{m\over 2}} \leq 1$. The inequality $\dim S(\mu_m) > 1$ is therefore equivalent to the fact that $\# \cD(m,a) \geq 2$. Choose two distinct pairs $\pm n$ and~$\pm n^\prime$ in $\cD(m,a)$. Note that if $m$ is even, we may have $n=-n$ or $n^\prime=-n^\prime$.

\subsection
\label{wpmncalc}
Choose $p$ and~$U$ as in Lemma~\ref{preW2lemma}, where we may further assume that the curves $C_t$ for $t \in U(k)$ are either all non-hyperelliptic or all hyperelliptic. We keep the notation introduced in the previous section; in particular we recall that $\cK(C_U/U) = \oplus_{n \in \mZ/m\mZ}\, \cK_{(n)}$.

It follows from~(\ref{-nai=m-nai}) and the definition of the forms~$\omega_{n,\nu}$ in~(\ref{omegannu}) that
$$
\eta := \omega_{n,0} \otimes \omega_{-n,0} - \omega_{n^\prime,0} \otimes \omega_{-n^\prime,0}
$$
is a section of~$\cK_{(0)}$. For $\nu \in \{\pm n,\pm n^\prime\}$ the matrix $A_\nu = A_\nu(t)$ of Lemma~\ref{InvgammaMatrix} is a polynomial in $\mF_p[t_1,\ldots,t_4]$. As the Jacobians of the curves~$C_t$ for $t\in U(k)$ are ordinary, the Cartier operator~$\gamma^{-1}$ in Lemma~\ref{InvgammaMatrix} is an isomorphism of vector bundles, so $A_\nu(t)$ is invertible as a section of~$\cO_U$. Because $\omega_{n,0} \cdot \omega_{-n,0} = \omega_{n^\prime,0} \cdot \omega_{-n^\prime,0}$ is a non-zero section of $f_*(\omega^{\otimes 2})$, it follows from Prop.~\ref{Criterion} that
$$
A_n \cdot A_{-n} = A_{n^\prime} \cdot A_{-n^\prime} \eqlabel{cc=cc}
$$
as polynomials.

Define $B_\nu := A_\nu|_{t_1=0}$, the polynomial obtained from~$A_\nu$ by substituting $t_1=0$.  Explicitly,
$$
B_\nu = (-1)^{p-1} \cdot \sum_{j_2+j_3+j_4=p-1}\, {q \cdot \repr{-\nu a_2}_m \choose j_2} \cdots {q \cdot \repr{-\nu a_4}_m \choose j_4} \cdot t_2^{j_2} t_3^{j_3} t_4^{j_4}\, .
\eqlabel{Bnformula}
$$
Fix an index $\iota \in \{2,3,4\}$, and write $\{1,2,3,4\} = \{1,\iota\} \amalg \{\kappa,\lambda\}$. Let $v_n(\iota)$ be the largest integer~$v$ such that $B_n$ is divisible by~$t_\iota^v$. From~(\ref{Bnformula}) and~(\ref{-nai=m-nai}) we find
$$
\eqalign{v_n(\iota) &= \max\bigl\{0,(p-1) - q\cdot \repr{-na_\kappa}_m - q\cdot \repr{-na_\lambda}_m\bigr\} \cr
&= \max\bigl\{0,q \cdot \repr{na_\kappa}_m + q \cdot \repr{n a_\lambda}_m - (p-1)\bigr\}\, ,\cr}
$$
and similarly for $v_{-n}(\iota)$.

Next let $w_{\pm n}(\iota)$ be the largest integer~$w$ such that $B_n \cdot B_{-n}$ is divisible by~$t_\iota^w$. We find
$$
\eqalign{w_{\pm n}(\iota) &= \max\bigl\{0,q \cdot \repr{na_\kappa}_m + q \cdot \repr{n a_\lambda}_m - (p-1)\bigr\}\cr
&\qquad + \max\bigl\{0,q \cdot \repr{-na_\kappa}_m + q \cdot \repr{-n a_\lambda}_m - (p-1)\bigr\}\, ,}
$$
which by~(\ref{-nai=m-nai}) we can rewrite as
$$
\eqalign{
w_{\pm n}(\iota) &= \max\bigl\{0,q \cdot \repr{na_\kappa}_m + q \cdot \repr{n a_\lambda}_m - (p-1)\bigr\}\cr
&\qquad + \max\bigl\{0,q \cdot \repr{na_1}_m + q \cdot \repr{n a_\iota}_m - (p-1)\bigr\}\cr
&= q \cdot \max\bigl\{\repr{na_1}_m + \repr{n a_\iota}_m, \repr{na_\kappa}_m + \repr{n a_\lambda}_m \bigr\} - (p-1)\, .\cr}\eqlabel{wformula}
$$ 

\subsection
{\it Lemma. --- Consider two $4$-tuples $(b_1,b_2,b_3,b_4)$ and $(b_1^\prime,b_2^\prime,b_3^\prime,b_4^\prime)$ in $\{1,\ldots,m-1\}^4$ with $\sum_{i=1}^4\, b_i = 2m = \sum_{i=1}^4\, b^\prime_i$. Suppose that for all partitions $\{1,2,3,4\} = \{1,\iota\} \amalg \{\kappa,\lambda\}$ we have $\{b_1+b_\iota,b_\kappa+b_\lambda\} = \{b^\prime_1+b^\prime_\iota,b^\prime_\kappa+b^\prime_\lambda\}$ as sets. Then there is an even permutation $\sigma \in A_4$ of order~$2$ such that either $b^\prime_i = b_{\sigma(i)}$ for all~$i$ or $b^\prime_i = m-b_{\sigma(i)}$ for all~$i$.}
\medskip

\Proof
Straightforward checking of all $8=2^3$ possible combinations (two for each partition). 
\QED

\subsection
The relation~(\ref{cc=cc}) implies that for any partition $\{1,2,3,4\} = \{1,\iota\} \amalg \{\kappa,\lambda\}$ we have $w_{\pm n}(\iota) = w_{\pm n^\prime}(\iota)$. Since $\sum_{i=1}^4\,  \repr{n a_i}_m = 2m = \sum_{i=1}^4\,  \repr{n^\prime a_i}_m$, it follows from~(\ref{wformula}) that
$$
\bigl\{\repr{na_1}_m + \repr{n a_\iota}_m, \repr{na_\kappa}_m + \repr{n a_\lambda}_m \bigr\} = \bigl\{\repr{n^\prime a_1}_m + \repr{n^\prime a_\iota}_m, \repr{n^\prime a_\kappa}_m + \repr{n^\prime a_\lambda}_m \bigr\}
$$
as sets. By the lemma it follows, possibly after replacing~$n$ by~$-n$, which we may do, that there is an even permutation $\sigma \in A_4$ of order~$2$ such that $\repr{n^\prime a_i}_m = \repr{n a_{\sigma(i)}}_m$ for all~$i$. As the pairs $\pm n$ and $\pm n^\prime$ are distinct, $\sigma \neq 1$. Possibly after a permutation of the indices we may therefore assume that
$$
\eqalign{
\repr{n^\prime a_1}_m = \repr{n a_2}_m\, ,\qquad & \repr{n^\prime a_3}_m = \repr{n a_4}_m\, ,\cr
\repr{n^\prime a_2}_m = \repr{n a_1}_m\, ,\qquad & \repr{n^\prime a_4}_m = \repr{n a_3}_m\, .\cr}\eqlabel{nai=naj}
$$
Since $\gcd(m,a_1,\ldots,a_4)=1$ this implies that $\gcd(m,n) = \gcd(m,n^\prime)$.

\subsection
Let $r = \gcd(n,m)$. Let $m^\prime = m/r$, and consider the triple $(m^\prime,N^\prime,a^\prime)$ obtained from $(m,N,a)$ by reduction modulo~$m^\prime$. Then we have a bijection $\cD(m,a) \isomarrow \cD(m^\prime,a^\prime)$ by $(\pm n) \mapsto (\pm {n\over r})$. Hence $\# \cD(m^\prime,a^\prime) \geq 2$. In particular, $N^\prime =4$, and $(m^\prime,N^\prime,a^\prime)$ is not equivalent to one of the triples listed in Table~\ref{tablespecials}. By our minimality assumption on~$(m,N,a)$, see~\ref{assumpt}, it therefore follows that $r=1$. Hence we may replace the $4$-tuple~$a$ by $n^\prime a = (n^\prime a_1,\ldots,n^\prime a_4)$, which means we may set $n^\prime =1$ in the above. In this case, (\ref{nai=naj}) tells us that $n^2 = 1$ in $(\mZ/m\mZ)$, and the $4$-tuple~$a$ has the form $a = (a_1,na_1,a_3,na_3)$ with $n \neq \pm 1$ and 
$$
(n+1)(a_1+a_3) = 2m\, .\eqlabel{=2mrel}
$$
{}From these relations we want to deduce that there is a pair $\pm\nu \in \cD(m,a)$ with $\gcd(m,\nu) \neq 1$, thereby obtaining a contradiction. Write $m = 2^k \cdot M$ with $M$ odd.

First assume $n \not\equiv -1$ modulo~$M$. For an odd prime number~$\ell$ and an exponent $e \geq 1$, the congruence $n^2 \equiv 1 \bmod{\ell^e}$ only has $n \equiv \pm 1$ as solutions. Hence there is an odd prime divisor~$\ell$ of~$m$ such that $n \not\equiv -1 \bmod{\ell}$. But then (\ref{=2mrel}) gives $a_1 + a_3 \equiv 0 \bmod{\ell}$, and since $\gcd(m,a_1,\ldots,a_4) = 1$ we further have $a_1 \not\equiv 0 \bmod \ell$ and $a_3 \not\equiv 0 \bmod \ell$. Taking $\nu := m/\ell$ we find that $\pm\nu \in \cD(m,a)$; contradiction.

The only remaining possibility is that $n \equiv -1$ modulo~$M$, and therefore $n \not\equiv -1 \bmod{2^k}$. In particular, $k \geq 2$. In this case (\ref{=2mrel}) implies $a_1 + a_3 \equiv 0 \bmod 4$, in which case $\pm m/4 \in \cD(m,a)$; contradiction. This completes the proof of Thm.~\ref{mainthm} in the case $N=4$.

\section{Proof of the main result: five or more branch points}{PfN>4}
\bigskip

\noindent
We now turn to the case $N\geq 5$. The idea behind the proof is the same as in the previous section, but there are some additional difficulties. In \ref{newpartnonrig} and~\ref{twopairs} we prove two technical results. Once we have these, we can again write down elements in the source of the map~(\ref{composite(0)}). Similar to what we did in Section~\ref{PfN=4}, this gives us polynomial identities~(\ref{cinvAcinvA}), and with some elementary combinatorial arguments we can then conclude the proof.

\subsection
\label{assumptN>4}
{\it Assumptions.\/} --- We retain the assumptions made in~\ref{assumpt}. In addition we now assume $N \geq 5$. Again we assume $(m,N,a)$ is minimal. Because we have already proven Thm.~\ref{mainthm} for $N=4$, it is again true that if $m^\prime$ is a proper divisor of~$m$ and $(m^\prime,N^\prime,a^\prime)$ is obtained from $(m,N,a)$ by reduction modulo~$m^\prime$, either $N^\prime \leq 3$ or $(m^\prime,N^\prime,a^\prime)$ is equivalent to one of the twenty triples listed in Table~\ref{tablespecials}.

\subsection
\label{newpartnonrig}
{\it Lemma. --- With assumptions as in\/~{\rm \ref{assumptN>4}}, there exists an index $n \in (\mZ/m\mZ)^*$ such that $\{d_n,d_{-n}\} \neq \{0,N-2\}$.}
\medskip

With the terminology introduced in~\ref{oldfactors}, the lemma says that the new part of the local system $\mV$ (as in Section~\ref{DecompVHS}) is not isotrivial.
\medskip

\Proof
We assume $\{d_n,d_{-n}\} = \{0,N-2\}$ for all $n \in (\mZ/m\mZ)^*$, and we seek to derive a contradiction. It will be convenient to consider the function $\delta_a \colon \mZ/m\mZ \to \mZ_{\geq 0}$ given by $\delta_a(n) = \sum_{i=1}^N\, \frpart{{-n a_i\over m}}$. Note that for $n \neq 0$ the $d_n$ defined in~(\ref{dndef}) equals $-1 + \delta_a(n)$. For $n \in (\mZ/m\mZ)^*$ we have $\delta_a(n) \geq 1$ and $\delta_a(n) + \delta_a(-n) = N$. Our assumption is equivalent to the condition that 
$$
\{\delta_a(n),\delta_a(-n)\} = \{1,N-1\}\quad \hbox{for all $n \in (\mZ/m\mZ)^*$.}
\eqlabel{deltandelta-n}
$$
Possibly after replacing $a = (a_1,\ldots,a_N)$ by~$-a$, we may assume $\delta_a(1) = 1$, which means that 
$$
\repr{a_1}_m + \cdots + \repr{a_N}_m = m\, .
\eqlabel{Sumai=m}
$$ 

We divide the proof into a couple of steps.
\medskip

\noindent
{\it Step 1.\/} Our first goal is to show that $m > 60$. For a given~$m$, (\ref{Sumai=m}) leaves only finitely many possible triples $(m,N,a)$, up to equivalence. With the help of a small computer program we check that for $m\leq 60$ the only two possibilities with $N \geq 5$ and $\{\delta_a(n),\delta_a(-n)\} = \{1,N-1\}$ for all $n \in (\mZ/m\mZ)^*$ are given (up to equivalence) by
$$
m=6\, ,\quad N=5\, ,\quad a = (1,1,1,1,2)\, ;
\eqlabel{65aExa}
$$
\vskip-2\bigskipamount
$$
m=6\, ,\quad N=6\, ,\quad a = (1,1,1,1,1,1)\, .
\eqlabel{66aExa}
$$
So it remains to show that for these triples, $Z(m,N,a)$ is not special.

In case (\ref{65aExa}), we have $g=7$ and the family of Jacobians $J \to T$ decomposes, up to isogeny, as a product of three factors: $J \sim Y_1 \times Y_2 \times Y_3$, with\par
\itemize
\item{---} $Y_1$ a pull-back of the Legendre family of elliptic curves, which is Example~{\oldstyle 1} in Table~\ref{tablespecials}; this is the old factor corresponding to the divisor $2$ of~$m$,
\item{---} $Y_2$ a pull-back of Example~{\oldstyle 6}; this is the old factor corresponding to the divisor $3$ of~$m$, 
\item{---} $Y_3$ the new part.
\stopitemize
\smallskip

\noindent
(The new part~$Y_3$ is in fact an isotrivial family of $3$-dimensional abelian varieties of CM type, but we will not need this.) Let $b \in T(\mC)$ be a Hodge-generic point. The fibre $Y_{1,b}$ is an elliptic curve with endomorphism ring~$\mZ$; its Hodge group is isomorphic with $\SL_2$. The fibre~$Y_{2,b}$ is an abelian threefold that is easily seen to be simple. (E.g., use the fact that Example~{\oldstyle 6} is special, together with [\ref{ShimuraAnFam}], Thm.~5.) By [\ref{BMYuZ}], Prop.~3.8 the Hodge group of $Y_{1,b} \times Y_{2,b}$ is the product of the two Hodge groups. Hence the smallest special subvariety $S \subset \cA_7$ that contains $Z(m,N,a)$ has dimension at least $1+2 = 3$ (see~\ref{MTfactors}), and since $N-3 =2 < 3$ this means that $Z(m,N,a)$ is not special.

In case (\ref{66aExa}), we have $g=10$ and the family of Jacobians $J \to T$ decomposes, up to isogeny, as a product of three factors: $J \sim Y_1 \times Y_2 \times Y_3$, with\par
\itemize
\item{---} $Y_1$ a pull-back of Example~{\oldstyle 2}; this is the old factor corresponding to the divisor $2$ of~$m$,
\item{---} $Y_2$ a pull-back of Example~{\oldstyle 10}; this is the old factor corresponding to the divisor $3$ of~$m$, 
\item{---} $Y_3$ the new part.
\stopitemize
\smallskip

\noindent
(In this case the new part~$Y_3$ is an isotrivial family of $4$-dimensional abelian varieties of CM type; again we will not need this.) For a Hodge-generic point $b \in T(\mC)$ the fibre $Y_{1,b}$ is an abelian surface with endomorphism ring~$\mZ$; its Hodge group $H_1 = \Hdg(Y_{1,b})$ is isomorphic with $\Sp_4$. The fibre~$Y_{2,b}$ is an abelian fourfold. By Thm.~1.1 of~[\ref{Zarhin}] we have $\End(Y_{2,b}) = \mZ[\zeta_3]$, and the tangent multiplicities for the factor~$Y_2$ (=Example~{\oldstyle 10}) are given by $d_{(1 \bmod 3)} = 3$ and $d_{(2\bmod 3)} = 1$. By [\ref{BMYuZDuke}], 7.4, it follows that for the generic Hodge group $H_2 = \Hdg(Y_{2,b})$ we have $H_{2,\mR} \cong \UU(1,3)$. If $H = \Hdg(J_b)$ is the generic Hodge group in our family of Jacobians, $H$ projects surjectively to~$H_1$ and~$H_2$. It follows that among the simple factors of~$H_\mR^\ad$ both a factor $\PSp_4$ and a factor~$\PSU(1,3)$ occur; this implies that the smallest special subvariety $S \subset \cA_7$ that contains $Z(m,N,a)$ has dimension at least $3+3 = 6$ (again see~\ref{MTfactors}), and since $N-3 =3 < 6$ this means that $Z(m,N,a)$ is not special. 
\medskip

\noindent
{\it Step 2.\/} The function $\delta_a$ has the property that $\delta_a(n_1 + n_2) \leq \delta_a(n_1) + \delta_a(n_2)$. In particular, if $n_1$, $n_2$ and $n_1 + n_2$ are all in $(\mZ/m\mZ)^*$ then by~(\ref{deltandelta-n}) we have the implication
$$
\delta_a(n_1) = \delta_a(n_2) = 1 \quad\Rightarrow\quad \delta_a(n_1+n_2) = 1\, .
\eqlabel{convexity}
$$
\medskip

\noindent
{\it Step 3.\/} Next we show that $m$ is not divisible by $2$, $3$ or~$5$. Let $p \in \{2,3,5\}$, and suppose $m$ is divisible by~$p$. Let $m^\prime = m/p$ and let $(m^\prime,N^\prime,a^\prime)$ be obtained from $(m,N,a)$ by reduction modulo~$m^\prime$. The new triple cannot be one of the examples listed in Table~\ref{tablespecials}, for this would imply (by inspection of the table) that $m = p \cdot m^\prime \leq 5 \cdot 12 = 60$, which we have excluded. Hence our minimality assumption on the triple $(m,N,a)$ implies that $N^\prime \leq 3$.

Since $N \geq 5$, there are at least two indices~$i$ such that $a_i \equiv 0$ modulo~$m^\prime$. For $p=2$ this contradicts~(\ref{Sumai=m}). Hence $2 \in (\mZ/m\mZ)^*$ and by repeated application of~(\ref{convexity}), starting from $\delta_a(1)=1$, we find that $\delta_a(2^k) =1$ for all~$k$.

For $p=3$ the only possibility left, in view of~(\ref{Sumai=m}), is that there are precisely two indices~$i$ with $\repr{a_i}_m = m/3$ and that $\repr{a_i}_m < m/3$ for all remaining indices. This contradicts the fact that $\delta_a(2) =1$. Hence $3$ does not divide~$m$.

Finally take $p=5$. The fact that $\delta_a(2) = 1$ implies that there is a unique index~$i$ with $\repr{a_i}_m > m/2$ and that $\repr{a_i}_m < m/2$ for all remaining indices. Since $N \geq 5$ and $N^\prime \leq 3$ this leaves (possibly after a permutation of the indices) two possibilities: either 
$$
\repr{a_1}_m = 3m/5\, ,\quad \repr{a_2}_m = m/5\, ,\quad \repr{a_i}_m < m/5\quad \hbox{for all $i >2$,}
$$ 
\vskip-\medskipamount\noindent
or 
$$
m/2 < \repr{a_1}_m < 3m/5\, ,\quad \repr{a_2}_m = \repr{a_3}_m = m/5\, ,\quad \repr{a_i}_m < m/5\quad \hbox{for all $i >3$.}
$$
In both cases we get a contradiction with the fact that $\delta_a(4) =1$.
\medskip

\noindent
{\it Step 4.\/} We now show that $m$ has at most four distinct prime divisors. To see this, let $p$ be any prime divisor of~$m$, let $m^\prime = m/p$ and let $(m^\prime,N^\prime,a^\prime)$ be obtained from $(m,N,a)$ by reduction modulo~$m^\prime$. If the new triple is one of the examples listed in Table~\ref{tablespecials}, it can only be Example~{\oldstyle 17}, for otherwise inspection of the table gives that $m^\prime$ (and hence~$m$) is divisible by $2$, $3$, or~$5$, which we have excluded. Example~{\oldstyle 17} has $N^\prime = 4$ branch points, so this fact, together with our minimality assumption on the triple $(m,N,a)$, implies that $N^\prime \leq 4$. 

The previous argument shows that for any prime divisor~$p$ of~$m$, there are at most~$4$ indices~$i$ such that $a_i$ is not divisible by $m/p$. On the other hand, $a_i \not\equiv 0 \bmod m$, so for a given index~$i$ there can be at most one prime divisor~$p$ of~$m$ such that $a_i$ is divisible by $m/p$. As $N \geq 5$ it follows that $m$ has at most four distinct prime divisors.
\medskip

\noindent
{\it Step 5; conclusion of the proof.\/} Combining the conclusion of Step~3 with~(\ref{convexity}), we find that for all $\rho \in \{1,2,3,4,5\}$ we have $(\rho \bmod m) \in (\mZ/m\mZ)^*$ and $\delta_a(\rho\bmod m) =1$.

Let $\nu$ be an integer such that $\gcd(m,\nu) = 1$ and $\delta_a(\nu \bmod m) =1$. Let $\rho$ be the smallest positive integer such that $\gcd(m,\nu+\rho) =1$. The fact that $m$ has at most four distinct prime divisors, all greater than~$5$, implies that $\rho \leq 5$. But then $\gcd(\rho,m) =1$ and from (\ref{convexity}) we get that $\delta_a(\nu+\rho \bmod m) = 1$. Repeating this, we find that $\delta_a(n) =1$ for all $n \in (\mZ/m\mZ)^*$, which contradicts the fact that $\delta_a(n) + \delta_a(-n) = N > 3$. This completes the proof of the lemma.
\QED

\subsection
\label{twopairs}
Our next goal is to show there exist two distinct pairs $\pm n$ and $\pm n^\prime$ in~$I(m)$ with
$$
d_{-n} = d_{-n^\prime} = 1\, ,\quad\hbox{and}\quad d_n = d_{n^\prime} = N-3\, .
\eqlabel{d-ndn2x}
$$
(In particular, $n\not\equiv 0$ and $n^\prime \not\equiv 0$.) For this we work over~$\mC$ and we consider the generic Mumford-Tate group~$M$ of the family $C/T$. As in~\ref{MTfactors}, let $M_\mR^\ad = Q_1 \times \cdots \times Q_r$ be the decomposition of~$M_\mR^\ad$ as a product of simple factors. The assumption that $Z(m,N,a)$ is special means that $N-3 = \sum_{i=1}^r\, d(Q_i)$, with notation as in~\ref{MTfactors}. Further, this assumption implies that the connected algebraic monodromy group $\Mono^0$ of the family is a normal subgroup of~$M$; in particular, $\Mono^{0,\ad}_\mR = \prod_{i\in K}\, Q_i$ for some subset $K \subset \{1,\ldots,r\}$.

By Lemma~\ref{newpartnonrig} there exists an index~$n$ in $(\mZ/m\mZ)^*$ with $\{d_n,d_{-n}\} \neq \{0,N-2\}$. By Prop.~\ref{monosurjects} it follows that one of the simple factors~$Q_i$, say~$Q_1$, is a $\PSU(d_n,d_{-n})$, which gives $d(Q_1) = d_n d_{-n}$. Because $n \in (\mZ/m\mZ)^*$ we have $d_n + d_{-n} = N-2$; hence $d_n d_{-n} \geq N-3$ with equality of and only if $\{d_n,d_{-n}\} = \{1,N-3\}$. Possibly after replacing $n$ by~$-n$ we therefore have $d_n = N-3$ and $d_{-n} =1$. Further, the relation $N-3 = \sum_{i=1}^r\, d(Q_i)$ implies that all other simple factors $Q_2,\ldots,Q_r$ are anisotropic (i.e., compact).

Our assumptions imply that there is another index pair $\pm n^\prime \in I(m)$ with $d_{n^\prime} \neq 0$ and $d_{-n^\prime} \neq 0$, for if this is not the case then $(m,N,a)$ satisfies (\ref{RohdeCond}), which implies it is one of the examples of Table~\ref{tablespecials}. If $m=2m^\prime$ is even and $n^\prime = m^\prime$ then by Prop.~\ref{monosurjects} one of the factors~$Q_i$ is a $\PSp_{2h}$ with $h = d_{n^\prime} > 0$. As $N-3 \geq 2$, this factor cannot be the factor $Q_1 = \PSU(1,N-3)$, and since $\PSp_{2h}$ is not anisotropic we arrive at a contradiction. So $n^\prime \neq -n^\prime$, and, again by Prop.~\ref{monosurjects}, one of the~$Q_i$ is a non-compact unitary factor $\PSU(d_{n^\prime},d_{-n^\prime})$. This factor must be the factor~$Q_1$; hence $\{d_{n^\prime},d_{-n^\prime}\} = \{1,N-3\}$, which gives what we want.

\subsection
Consider two distinct pairs $\pm n$ and $\pm n^\prime$ in $I(m)$ for which~(\ref{d-ndn2x}) holds. Let $\Gamma \in \GL_{N-3}(\cO_U)$ be the matrix of $\gamma_{(n)} \colon F_U^* \mE_{U,(n)} \to \mE_{U,(n)}$ with regard to the frames given by the forms $\omega_{n,\sigma}$. It is the inverse of the matrix $A = A_n$ of Lemma~\ref{InvgammaMatrix}. Similarly, let $\Gamma^\prime$ be the matrix of~$\gamma_{(n^\prime)}$, and let $A^\prime = A_{n^\prime}$ be its inverse. Further, let $c$ and~$c^\prime$ be the sections of~$\cO_U^*$ (the inverses of the $1\times 1$ matrices $A_{-n}$ and~$A_{-n^\prime}$ of Lemma~\ref{InvgammaMatrix}) such that $\gamma(\omega_{-n,0} \otimes 1) = c \cdot \omega_{-n,0}$ and $\gamma(\omega_{-n^\prime,0} \otimes 1) = c \cdot \omega_{-n^\prime,0}$.

For $\sigma \in \{0,1,\ldots,N-4\}$, let $\eta_\sigma \in \Gamma(U,\cK_{(0)})$ be given by 
$$
\eta_\sigma := \omega_{-n,0} \otimes \omega_{n,\sigma} - \omega_{-n^\prime,0} \otimes \omega_{n^\prime,\sigma}\, .
$$
As $\omega_{-n,0} \cdot \omega_{n,\rho} = \omega_{-n^\prime,0} \cdot \omega_{n^\prime,\rho}$ as sections of $f_*(\omega^{\otimes 2})$, the image of~$\eta_\sigma$ under~$\Sym^2(\gamma)$ equals
$$
\sum_{\rho=0}^{N-4} \left(c \cdot \Gamma_{\rho,\sigma} - c^\prime \cdot \Gamma^\prime_{\rho,\sigma}\right) \cdot (\omega_{-n,0} \cdot \omega_{n,\rho})\, .
$$
Because the sections $\omega_{-n,0} \cdot \omega_{n,\rho}$ for $\rho \in \{0,\ldots,N-4\}$ are linearly independent, it follows from Prop.~\ref{Criterion} that $c \cdot \Gamma_{\rho,\sigma} - c^\prime \cdot \Gamma^\prime_{\rho,\sigma}$ for all $\sigma$, $\rho \in \{0,1,\ldots,N-4\}$, which can be rewritten as
$$
c^{-1} \cdot A_{\rho,\sigma} = c^{\prime,-1} \cdot A^\prime_{\rho,\sigma}\, .
\eqlabel{cinvAcinvA}
$$

Choose two distinct indices $\kappa$, $\lambda$ in $\{1,2,\ldots,N\}$, and write $\{1,\ldots,N\} = \{\kappa,\lambda\} \amalg I$. For $\nu \in \{0,1\}$, define $v_n(\nu)$ to be the largest integer~$v$ such that $A_{\nu,\nu}|_{t_\kappa=0}$ is divisible by~$t_\lambda^v$. Similarly, let $v_{-n}$ be the largest integer~$v$ such that $c^{-1}|_{t_\kappa=0}$ is divisible by~$t_\lambda^v$. Using the explicit formulas for $c^{-1}$ and the matrix~$A$ from Lemma~\ref{InvgammaMatrix}, we find
$$
\eqalign{v_{-n} &= \max\Bigl\{0,(p-1) - q \cdot\sum_{i\in I}\, \repr{na_i}_m \Bigr\}\, ,\cr
v_n(0) &= \max\Bigl\{0,(N-3)(p-1) - q\cdot \sum_{i\in I}\, \repr{-na_i}_m \Bigr\}\, ,\cr
v_n(1) &= 0\, .\cr}
$$
(Recall that $q = (p-1)/m$.)

Next we define $w_n(\nu) = v_{-n} + v_n(\nu)$ to be the largest integer~$w$ such that $\bigl(c^{-1} \cdot A_{\nu,\nu}\bigr)|_{t_\kappa=0}$ is divisible by~$t_\lambda^w$. Similar to~(\ref{-nai=m-nai}), it follows from the fact that $d_{-n} =1$ and $d_n = N-3$ that $\repr{-na_i}_m = m - \repr{na_i}_m$ for all~$i$ and $\sum_{i=1}^N\, \repr{n a_i}_m = 2m$. Using these relations we obtain
$$
\eqalign{
w(0) &= \max\Bigl\{(p-1) - q \cdot \sum_{i \notin I}\, \repr{n a_i}_m, (p-1) - q \cdot \sum_{i \in I}\, \repr{n a_i}_m \Bigr\}\cr
&= \Bigl\vert (p-1) - q \cdot\sum_{i\notin I}\, \repr{na_i}_m \Bigr\vert = q \cdot \bigl\vert m - \repr{na_\kappa}_m - \repr{na_\lambda}_m\bigr\vert\, ,\cr
w(1) &= \max\Bigl\{0,(p-1) - q \cdot\sum_{i\in I}\, \repr{na_i}_m \Bigr\}\cr
&= \max\Bigl\{0,(p-1) - q \cdot\sum_{i\notin I}\, \repr{na_i}_m \Bigr\} = q\cdot  \max\bigl\{0,\repr{na_\kappa}_m + \repr{na_\lambda}_m -m \bigr\}\, .\cr}
$$

We can do the same calculations for the pair $\pm n^\prime$; let us call $w^\prime(0)$ and~$w^\prime(1)$ the resulting values. {}From~(\ref{cinvAcinvA}) we get the relations $w(0) = w^\prime(0)$ and $w(1) = w^\prime(1)$. It follows that for all choices of $\kappa$ and~$\lambda$ we have $\repr{n a_\kappa}_m + \repr{na_\lambda}_m = \repr{n^\prime a_\kappa}_m + \repr{n^\prime a_\lambda}_m$. This readily implies that $\repr{n a_i}_m = \repr{n^\prime a_i}_m$ for all $i \in \{1,\ldots,N\}$. As $\gcd(m,a_1,\ldots,a_N) = 1$ it follows that $n=n^\prime$, which contradicts our assumption that the pairs $\pm n$ and~$\pm n^\prime$ are distinct. This completes the proof of Theorem~\ref{mainthm}.
\QED

\vskip2.0\bigskipamount plus 2pt minus 1pt%
\goodbreak\centerline{{\sectitlefont References}}%
\nobreak\vskip.75\bigskipamount plus 2pt minus 1pt%

\Reference{ACampo}
N.~A'Campo, {\it Tresses, monodromie et le groupe symplectique.} Comment.\ Math.\ Helv.\ 54 (1979), 318--327.

\Reference{AchtPries}
J.~Achter, R.~Pries, {\it The integral monodromy of hyperelliptic and trielliptic curves.} Math.\ Ann.\ 338 (2007), 187--206. 

\Reference{Andre}
Y.~Andr\'e, {\it Mumford-Tate groups of mixed Hodge structures and the theorem of the fixed part.} Compositio Math.\ 82 (1992), 1--24.

\Reference{Bouw}
I.~Bouw, {\it The $p$-rank of ramified covers of curves.} Compositio Math.\ 126 (2001), 295--322.

\Reference{ChevWeil}
C.~Chevalley, A.~Weil, {\it \"Uber das Verhalten der Integrale 1.\ Gattung bei 
Automorphismen des Funktio\-nen\-k\"orpers.} Abh.\ Math.\ Semin.\ Univ.\ Hambg.\ 10 (1934), 358--361.

\Reference{dJN}
J.~de Jong, R.~Noot, {\it Jacobians with complex multiplication.} In: Arithmetic algebraic geometry (Texel, 1989), Progr.\ Math.\ 89, Birkh\"auser, 1991, pp.\ 177--192.

\Reference{DelMost}
P.~Deligne, G.~Mostow, {\it Monodromy of hypergeometric functions and nonlattice integral monodromy.} Inst. Hautes \'Etudes Sci.\ Publ.\ Math.\ 63 (1986), 5--89. 

\Reference{DworkOgus}
B.~Dwork, A.~Ogus, {\it Canonical liftings of Jacobians.} Compositio Math.\ 58 (1986), 111--131. 

\Reference{Katz}
N.~Katz, {\it Serre-Tate local moduli.} In: Algebraic surfaces, Lecture Notes in Math.\ 868, Springer-Verlag, 1981, pp.\ 138--202.

\Reference{Messing}
W.~Messing, {\it The crystals associated to Barsotti-Tate groups: with applications to abelian schemes.\/} Lecture Notes in Math.\ 264, Springer-Verlag, 1972.

\Reference{BMthesis}
B.~Moonen, {\it Special points and linearity properties of Shimura varieties.} PhD thesis, University of Utrecht, 1995.

\Reference{BMLin1}
B.~Moonen, {\it Linearity properties of Shimura varieties.~I.}
J.\ Alg.\ Geom.\ 7 (1998), 539--567. 

\Reference{BMLin2}
B.~Moonen, {\it Linearity properties of Shimura varieties.~II.}
Compositio Math.\ 114 (1998), 3--35. 

\Reference{BMFO}
B.~Moonen, F.~Oort, {\it The Torelli locus and special subvarieties.} Preprint, August 2010, 43pp.

\Reference{BMYuZDuke}
B.~Moonen, Yu.~Zarhin, {\it Hodge classes and Tate classes on simple abelian fourfolds.} Duke Math.\ J. 77 (1995), 553--581. 

\Reference{BMYuZ}
B.~Moonen, Yu.~Zarhin, {\it Hodge classes on abelian varieties of low dimension.} Math.\ Ann.\ 315 (1999), 711--733. 

\Reference{Most}
G.~Mostow, {\it On discontinuous action of monodromy groups on the complex $n$-ball.} J.\ Amer.\ Math.\ Soc.\ 1 (1988), 555--586.

\Reference{NootThesis}
R.~Noot, {\it Hodge classes, Tate classes, and local moduli of abelian varieties.} PhD thesis, University of Utrecht, 1992. 

\Reference{NootModels}
R.~Noot, {\it Models of Shimura varieties in mixed characteristic.}
J.\ Alg.\ Geom.\ 5 (1996), 187--207.

\Reference{NootBourb}
R.~Noot, {\it Correspondances de Hecke, action de Galois et la conjecture d'Andr\'e-Oort (d'apr\`es Edixhoven et Yafaev).} S\'em.\ Bourbaki, Exp.\ 942. Ast\'erisque 307 (2006), 165--197. 

\Reference{FOJS}
F.~Oort, J.~Steenbrink, {\it The local Torelli problem for algebraic curves.} In:  Journ\'ees de G\'eometrie Alg\'ebrique d'Angers, Sijthoff \& Noordhoff, 1980, pp.\ 157--204. 

\Reference{Rohde}
J.~Rohde, {\it Cyclic coverings, Calabi-Yau manifolds and complex multiplication.} Lecture Notes in Math.\ 1975, Springer, 2009.

\Reference{ShimuraAnFam}
G.~Shimura, {\it On analytic families of polarized abelian varieties and automorphic functions.} Ann. of Math. (2) 78 (1963), 149--192. 

\Reference{Shimura}
G.~Shimura, {\it On purely transcendental fields of automorphic functions of several variable.} Osaka J.\ Math.\ 1 (1964), 1--14. 

\Reference{Yafaev}
A.~Yafaev, {\it A conjecture of Yves Andr\'e's.} Duke Math.\ J. 132 (2006), 393--407.

\Reference{Zarhin}
Yu.~Zarhin, {\it The endomorphism rings of Jacobians of cyclic covers of the projective line.} Math.\ Proc.\ Cambridge Philos.\ Soc.\ 136 (2004), 257--267.

\Addresses
\end